\documentclass[10pt,sumlimits,namelimits]{article}

\usepackage{amsthm}
\usepackage{amsmath}
\usepackage{amsfonts}
\usepackage{amssymb}
\usepackage{color}
\usepackage[colorlinks=true,urlcolor=blue]{hyperref}
\usepackage[applemac]{inputenc}

\newcommand{\N}{\mathbb{N}}
\newcommand{\Z}{\mathbb{Z}}

\newcommand{\Q}{\mathbb{Q}}
\newcommand{\R}{\mathbb{R}}
\newcommand{\C}{\mathbb{C}}
\newcommand{\vA}{{\cal A}}
\newcommand{\vC}{{\cal C}}
\newcommand{\vD}{{\cal D}}

\newcommand{\vG}{{\cal G}}
\newcommand{\vH}{{\cal H}}
\newcommand{\vL}{{\cal L}}
\newcommand{\vN}{{\cal N}}
\newcommand{\vP}{{\cal P}}
\newcommand{\vS}{{\cal S}}
\newcommand{\vphi}{\varphi}
\newcommand{\eps}{\varepsilon}
\newcommand{\dsp}{\displaystyle}
\newcommand{\ovl}{\overline}
\newcommand{\udl}{\underline}
\newcommand{\vlim}{\lim\limits}

\newcommand{\vliminf}{\liminf\limits}

\newcommand{\vinf}{\inf\limits}
\newcommand{\vint}{\int\limits}
\newcommand{\vsum}{\sum\limits}

\newcommand{\tends}{\longrightarrow}

\newcommand{\loc}{{\rm loc}}
\newcommand{\Id}{{\rm Id}}
\renewcommand{\b}{{\rm b}}
\renewcommand{\r}{{\rm r}}
\renewcommand{\le}{\leqslant}
\renewcommand{\ge}{\geqslant}

\numberwithin{equation}{section}

\newtheorem{thm}{Theorem}[section]
\newtheorem{prop}[thm]{Proposition}

\newtheorem{lem}[thm]{Lemma}

\newtheorem{vthm}{Theorem}[subsection]
\newtheorem{vprop}[vthm]{Proposition}

\newtheorem{vlem}[vthm]{Lemma}

\theoremstyle{definition}
\newtheorem{rmk}[thm]{Remark}
\newtheorem{defi}[thm]{Definition}

\newtheorem{exa}[thm]{Example}

\theoremstyle{definition}
\newtheorem{vrmk}[vthm]{Remark}
\newtheorem{vdefi}[vthm]{Definition}
\newtheorem{vnota}[vthm]{Notation}

\newenvironment{proof*}{\noindent{\bf Proof.}}{\qed}
\newenvironment{vproof}[1]{\noindent{\bf Proof #1}}{\qed}

\title{\Huge \sc A Generalized Interpolation Inequality and its Application to the Stabilization of Damped Equations}
\author{\sc Pascal Bégout\footnote{Supported by Grant HPRN--CT--2002-00284 of the European program {\it New materials, adaptive systems and their nonlinearities: modelling, control and numerical simulation}} and Fernando Soria\footnote{Supported by Grant BFM2001--0189 of the MCYT (Spain)}}
\date{}

\begin{document}

\maketitle

$$
\begin{array}{cc}
^*\mbox{Laboratoire d'Analyse et Probabilit\'es } & \;^\dagger\mbox{Departamento de Matem\'aticas} \\
           \mbox{D\'epartement de Math\'ematiques } & \mbox{ Facultad de Ciencias} \\
            \mbox{Universit\'e d'\'Evry Val d'Essonne } & \mbox{ Universidad Aut\'onoma de Madrid} \\
                  \mbox{Boulevard François Mitterrand } & \mbox{ Cantoblanco} \\
                   \mbox{91025 \'Evry Cedex, FRANCE } & \mbox{ 28049 Madrid, SPAIN}
\bigskip \\
\mbox{
{\footnotesize$^*$e-mail\:: }\htmladdnormallink{{\footnotesize
\udl{\tt{Pascal.Begout@univ-evry.fr}}}}{mailto:Pascal.Begout@univ-evry.fr} } & \mbox{{\footnotesize$^\dagger$e-mail\:: }
\htmladdnormallink{{\footnotesize\udl{\tt{fernando.soria@uam.es}}}}
{mailto:fernando.soria@uam.es}}
\bigskip
\end{array}
$$

\begin{abstract}
In this paper, we establish a generalized Hölder's or interpolation inequality for weighted spaces in which the weights are non-necessarily homogeneous. We apply it to the stabilization of some damped wave-like evolution equations. This allows obtaining explicit decay rates for smooth solutions for more general classes of damping operators. In particular, for $1-d$ models, we can give an explicit decay estimate for pointwise damping mechanisms supported on any strategic point.
\end{abstract}

{\let\thefootnote\relax\footnotetext{2000 Mathematics Subject Classification: 93D15 (35B37, 35L70, 93B52, 93C20)
}}
{\let\thefootnote\relax\footnotetext{Key Words: damped equations, damping control, generalized Hölder's inequality, interpolation inequality, stabilization
}}

\tableofcontents

\baselineskip .7cm

\section{Introduction}
\label{intro}

We are interested on a generalized Hölder's or interpolation inequality, in order to establish explicit decay rates for smooth solutions of damped wave-like equations with weak damping.
\medskip
\\
Let $(\Omega,\Upsilon,\mu)$ be a measure space and let $\omega_1$ and $\omega_2$ be two $\mu$-measurable  weights on $\Omega.$ The problem we address consists in finding suitable functions $\Phi$ and $\Psi$ such that
\begin{gather}
 \label{intro1}
  1\le\Phi\left(\dfrac{\dsp\int_\Omega|f(x)|\omega_1(x)d\mu(x)}{\|f\|_{L^1(\Omega,\Upsilon,\mu)}}\right)
  \Psi\left(\dfrac{\dsp\int_\Omega|f(x)|\omega_2(x)d\mu(x)}{\|f\|_{L^1(\Omega,\Upsilon,\mu)}}\right),
\end{gather}
for any $f\in L^1(\Omega,\Upsilon,\mu)\cap L^1(\Omega,\Upsilon,\omega_1d\mu)\cap L^1(\Omega,\Upsilon,\omega_2d\mu).$
\medskip
\\
The case where the weights functions are homogeneous is well-known. Indeed, if $\omega_1(x)=|x|^\alpha$ and $\omega_2(x)=|x|^{-\beta}$ $(\alpha,\beta>0),$ the classical Hölder's inequality gives
\begin{gather}
 \label{intro2}
  \int_\Omega|f(x)|dx\le\left(\int_\Omega|f(x)|\:|x|^\alpha dx\right)^\frac{\beta}{\alpha+\beta}
  \left(\int_\Omega|f(x)|\:|x|^{-\beta}dx\right)^\frac{\alpha}{\alpha+\beta},
\end{gather}
where $dx$ denotes the Lebesgue's measure or, equivalently,
\begin{gather*}
 1\le\left(\dfrac{\dsp\int_\Omega|f(x)|\:|x|^\alpha dx}{\dsp\int_\Omega|f(x)|dx}\right)^\frac{\beta}{\alpha+\beta}
 \left(\dfrac{\dsp\int_\Omega|f(x)|\:|x|^{-\beta}dx}{\dsp\int_\Omega|f(x)|dx}\right)^\frac{\alpha}{\alpha+\beta}.
\end{gather*}

\noindent
Obviously, (\ref{intro2}) is a particular case of (\ref{intro1}), in which the functions $\Phi$ and $\Psi$ are respectively $\Phi(t)=t^\frac{\beta}{\alpha + \beta}$ and $\Psi(t)=t^\frac{\alpha}{\alpha + \beta}$.
\medskip
\\
This paper is devoted to obtain a generalization of (\ref{intro2}) for non-homogeneous weights. We are typically interested in situations in which, for instance,  $\omega_1(x)=e^{-|x|}$ and $\omega_2(x)=|x|^2$. As we shall see, if we are able to get an interpolation inequality of the form (\ref{intro1}) in this case, we will be able to give new explicit decay rates for damped $1-d$ wave equations with pointwise damping.
\medskip
\\
Let us briefly illustrate the connection between these two issues.
\medskip
\\
Let $a\in L^\infty(0,1)$ be a nonnegative and bounded damping potential and consider the damped wave equation in one space dimension,
\begin{equation}
 \label{intro3}
  \left\{
   \begin{array}{rcl}
     u_{tt}(t,x)-u_{xx}(t,x)+a(x)u_t(t,x)=0, & \mbox{for} & (t,x)\in(0,\infty)\times(0,1), \medskip \\
                                         u(t,0)=u(t,1)=0, & \mbox{for} & t\in[0,\infty), \medskip \\
          u(0,x)=u^0(x), \: u_t(0,x)=u^1(x), & \mbox{for} & x\in(0,1).
   \end{array}
  \right.
\end{equation}

\noindent
This system is well-posed. More precisely, for any initial data $u^0\in H^1_0(0,1)$ and $u^1\in L^2(0,1),$ there exists a unique solution in the class $\vC([0,\infty); H^1_0(0, 1))\cap\vC^1([0,\infty); L^2(0,1))$. The energy of solutions
\begin{gather*}
E(t)=\dfrac{1}{2}\left(\|u_t(t)\|_{L^2(0,1)}^2+\|u_x(t)\|_{L^2(0,1)}^2\right),
\end{gather*}
decreases along trajectories according to the dissipation law
\begin{gather}
\label{dr}
 \frac{d}{dt}E(t)=-\vint_0^1 a(x) |u_t(t,x)|^2 dx.
\end{gather}

\noindent
The decay rate of the energy depends on the efficiency of the damping term when absorbing the energy of the system according to (\ref{dr}).
\medskip
\\
Using LaSalle's invariance principle, it is easy to see that the energy of every solution tends to zero as $t\tends\infty$ whenever the damping potential $a$ satisfies for almost every $x\in I,$ $a(x)\ge a_0>0,$ for some constant $a_0>0,$ where $I\subset(0,1)$ is a set of positive measure (Haraux~\cite{MR804885}). In the $1-d$ case under consideration, in fact, one can even show that the energy of solutions tends to zero exponentially. To prove this fact, it is sufficient to show that for some $T >0$ and $C>0$ the following inequality holds
\begin{gather}
 \label{oi}
  E(0)\le C\vint_0^T\!\!\!\vint_0^1a(x) |u_t(t,x)|^2 dxdt,
\end{gather}
for every solution.
\medskip
\\
This inequality, which is often referred to as {\it observability inequality}, asserts that the damping mechanism during a time interval $(0, T)$ suffices to capture a fraction of the total energy of all solutions.
\medskip
\\
Combining (\ref{dr}), (\ref{oi}) and the semigroup property, it is easy to see that the exponential decay property holds, i.e. there exist $C>0$ and $\omega >0$ such that 
\begin{gather}
 \label{dp}
  \forall t \ge 0, \; E(t)\le C E(0) e^{-\omega t},
\end{gather}
for every solution.
\medskip
\\
In fact, to prove that (\ref{oi}) is fulfilled, one can use the fact that it is sufficient to prove it for the solutions of the corresponding conservative systems (\ref{intro3}) with $a=0.$ In that case, the inequality is easy to get for $T=2$ using the Fourier decomposition of solutions. 
\medskip
\\
Let us now consider a case where the control is supported simply on a  point $a\in(0,1)$ through a Dirac mass,
\begin{gather}
 \label{intro5}
  u_{tt}-u_{xx}+\delta_au_t(t,a)=0, \quad (t,x)\in(0,\infty)\times(0,1),
\end{gather}
with the same boundary conditions, initial data and energy as before. Here, $\delta_a$ denotes the Dirac mass concentrated in $a.$ 
\medskip
\\
When the point $a \in \Q$, there are solutions of (\ref{intro5}) that do not decay and for which the energy is constant in time. This is due to the fact that rational points are nodal ones for the corresponding Sturm-Liouville problem. 
\medskip
\\
When $a\not\in\Q,$ LaSalle's invariance principle allows proving that the energy of each solution tends to zero as $t\tends\infty$. However, in this case the exponential decay rate does not hold. This is due to the fact that, even if $a\not\in\Q,$ the damping term does not dissipative uniformly all the Fourier components of the solutions. This can be easily seen when analyzing the analogue of (\ref{oi}). Indeed,  there exists a sequence of separate variable solutions of the conservative problem (\ref{intro3}) with $a=0$ for which the energy $E(0)$ is of order one and the dissipated quantity, $\dsp\int_0^T |u_t(t,a)|^2 dt,$ tends to zero. This sequence can be built in separated variables, based on the sequence of eigenfunction $\sin(nx)$ such that $\sin (n a)$ tends to zero as $n$ tends to infinity. The main difference with the case where the damping potential $a\ge0$ is positive on a set of positive measure is that, in that case, $\vinf_{n\ge 1} \dsp\int_0^1 a(x) \sin^2(nx) dx >0.$
\medskip
\\
In view of this, one may only expect a weaker observability inequality to hold. A natural way of proceeding in this case is to obtain a weakened version of (\ref{oi}) in which the  energy $E(0)$ in the left hand side is replaced  by a weaker energy $E_-(0)$ which, roughly speaking, is the Fourier norm of solutions with weights $\sin^2 (n a)$. More precisely, 
\begin{gather}
 \label{oiweak}
  E_-(0)\le C\vint_0^T |u_t(t,a)|^2 dt= - C(E(T)-E(0)).
\end{gather}

\noindent
The problem is then how to derive an explicit decay rate for the energy $E$ out of (\ref{oiweak}). First, we need to assume some more regularity on the initial data, say, $(u^0, u^1)\in [H^2(0, 1)\cap H^1_0(0,1)]\times H^1_0(0, 1).$  We denote by $E_+$ the corresponding energy, $E_+(0)=\frac{1}{2}\|(u^0,u^1)\|^2_{H^2(0,1)\times H^1_0(0,1)}.$
\medskip
\\
In this way, we have three different energies with different degrees of strength: $E$, which is the reference energy in which we are interested, $E_+,$ which is finite because the initial data have been taken to be smooth, and $E_-$ which is the weaker energy the damping really damps out according to (\ref{oiweak}).
\medskip
\\
Applying (\ref{intro1}), one can deduce an interpolation inequality of the form

\begin{gather}
 \label{intro7}
  1\le\Phi\left(\dfrac{E_-(0)}{E(0)}\right)\Psi\left(\dfrac{E_+(0)}{E(0)}\right),
\end{gather}
where $\Phi$ and $\Psi$ depend on the energies $E_+$ and $E_-$ under consideration, $E_+(0)$ being the strong norm $E_+(0)=\frac{1}{2}\|(u^0,u^1)\|^2_{H^2(0,1)\times H^1_0(0,1)}$.
This clearly implies
\begin{gather}
 \label{intro7b}
  E(0) \Phi^{-1} \left(\frac{1}{\Psi\left(\frac{E_+(0)}{E(0)}\right)}\right) \le E_-(0),
\end{gather}
which, together with the weak observability inequality (\ref{oiweak}) yields,
\begin{gather}
 \label{intro8}
  E(0)\Phi^{-1}\left(\frac{1}{\Psi\left(\frac{E_+(0)}{E(0)}\right)}\right) \le C(E(0)-E(T)),
\end{gather}
which, together with the semigroup property yield (see Ammari and Tucsnak~\cite{MR2002f:93104}),
\begin{gather}
 \label{drn}
  \forall t\ge0, \; E(t)\le
  \frac{C}{\Psi^{-1}\left(\frac{1}{\Phi\left(\frac{1}{t+1}\right)}\right)}\|(u^0,u^1)\|_{H^2(0,1)\times H^1_0(0,1)}^2.
\end{gather}

\noindent
Our method is closely of that one developed by Nicaise~\cite{MR2044992}, in which the decay estimate of the energy looks like \eqref{drn} (see Section~5 in~\cite{MR2044992}). But unfortunately, his method cannot apply in this paper because the damping term has to be more regular, in some sense, that one we consider (see \cite{MR2044992}).
\medskip
\\
Obviously, the decay rate in (\ref{drn}) depends on the behavior of the functions $\Psi$ and $\Phi$. More precisely, it depends on the behavior of $\Phi(t)$ near $t=0$ and then of that of $\Psi^{-1}$ at infinity. Therefore, in order to determine the decay of solutions it is necessary to have a sharp description of the functions $\Phi$ and $\Psi$ entering in the interpolation inequality. 
\medskip
\\
The behavior of $\Phi$ and $\Psi$ depends on the energies $E,$ $E_+$ and $E_-$ under consideration. We recall that $E_-$ is given by the weak observability inequality (\ref{oiweak}). This is intimately related to the weakness of the damping mechanism and no choice can be done at that level. By the contrary, there is some liberty at the level of choosing $E_+$  since the initial data can be chosen to be as smooth as we like. Obviously, one expects a faster decay rate for solutions when they are smoother. This is indeed the case as our analysis shows. All this can be precisely quantified by the analysis of the functions $\Phi$ and $\Psi$ in the interpolation inequality.
\medskip
\\
How $\Phi$ and $\Psi$ depend on the energies $E_+$ and $E_-,$ in the general context of the interpolation inequality (\ref{intro1}), corresponds to analyzing how the functions $\Phi$ and $\Psi$ depend on the weight functions $\omega_1$ and $\omega_2.$ This article is precisely devoted to prove a rather general version of (\ref{intro1}) with a careful analysis of the behavior of $\Phi$ and $\Psi.$ This will allow us to get explicit decay rates not only for the model problem above of the $1-d$ wave equation with pointwise damping but also for some other models that we shall discuss below. In particular, we will be able to give explicit decay rates for the stabilization of a beam by means of a piezoelectric actuators, a problem that was discussed by Tucsnak~\cite{MR95g:93039,MR97c:93024} in the context of control.
\medskip
\\
There is an extensive literature concerning the stabilization of damped wave-like equations. But most of it refers to the case where the damping term (linear or nonlinear one) is able to capture the whole energy of the system (see, for instance, Haraux and Zuazua~\cite{MR913963}, Nicaise~\cite{MR2044992} and Zuazua~\cite{MR92c:35021}). In these works, the multiplier method is implied, as a tool to quantify the amount of energy that the dissipative mechanism is able to observe. But to apply this method, the damping term has to be active in a large subset of the domain or of the boundary where the equation holds. Much less is known when the damping term is located in a narrow set, like, for instance, pointwise dampers in one space dimension. But, as we have shown above, the results one may expect in that setting need to be necessarily of a weaker nature since in those situations the damping term is only able to absorb a lower order energy. In particular, in this context, multiplier methods do not apply.
\medskip
\\
We focus mainly on the wave equation with a damping control concentrated on an interior point. Some partial results of explicit decay rates already exist and can be found in Ammari, Henrot and Tucsnak~\cite{MR2000m:93078,MR2002j:93073}, Jaffard, Tucsnak and Zuazua~\cite{MR99g:93073} and Tucsnak~\cite{MR1627724}. As explained above, our generalized interpolation inequality allows answering to this in much more generality. We will also address  the stabilization of Bernoulli--Euler beams with force and moment damping. For partial results of explicit decay rates, see Ammari and Tucsnak~\cite{MR1814271}. 
\medskip
\\
This paper is organized as follows. In Section~\ref{interpolation}, we establish our generalized Hölder's inequality or interpolation inequality (Theorems~\ref{ineg1} and \ref{thmineg1}). In Section~\ref{optimality}, we give a criterion of optimality for Theorem~\ref{ineg1} (Definition~\ref{defioptim}) and a sufficient condition to have optimality in our interpolation inequality (Proposition~\ref{propoptim1}). In Section~\ref{damping}, we apply these results to get explicit decay rates for the damped wave (see (\ref{wave})) with Dirichlet boundary condition and in Section~\ref{damping2} we briefly explain how these results can be applied to the wave equation with mixed boundary condition (Subsection~\ref{sswaves}, equation (\ref{waves})) and to some beam equations (Subsection~\ref{ssbeam}, equation (\ref{beam})). The explicit decay rates are given. These results extend the previous ones  by Ammari, Henrot and Tucsnak~\cite{MR2002j:93073}, Ammari and Tucsnak~\cite{MR1814271} and Jaffard, Tucsnak and Zuazua~\cite{MR99g:93073}.
\medskip
\\
We end this section by introducing some notations. For a real valued function $f$ defined on an open interval $I$ (respectively, $(m,\infty)$ for some $m\in\R)$ and for $a\in\partial I$ (respectively, $a\in\{m,\infty\}),$ the notation $f(a)$ means $\vlim_{\underset{t\in I}{t\to a}}f(t).$ For $a\in\R,$ we denote by $\delta_a$ the Dirac mass concentrated in $a.$

\section{An interpolation inequality}
\label{interpolation}

\noindent
Our analysis requires some elementary notions and results on convex functions. 
\medskip
\\
Recall that if $f:I\tends\R$ is a convex function on an open interval $I,$ then it is continuous,  locally absolutely continuous on $I$ and it is of class $\vC^1$ almost everywhere. More precisely, there exists a finite or countable set $\vN\subset I$ such that $f$ is of class $\vC^1$ relatively to $I\setminus\vN.$ In particular, for any $t,s\in I,$
$f(t)-f(s)=\dsp\int_s^tf^\prime(\sigma)d\sigma.$ In addition, $f^\prime$ is nondecreasing relatively to $I\setminus\vN.$ Furthermore, $f$ has a left derivative $f^\prime_\ell$ and a right derivative $f^\prime_\r$ at each point of $I$ and for any $t,s\in I$ such that $s<t,$ $f^\prime_\ell(s)\le f^\prime_\r(s)\le f^\prime_\ell(t)\le f^\prime_\r(t).$ For more details, see Niculescu and Persson~\cite{ni} (Theorems~1.3.1 and 1.3.3, p.12, Proposition~3.4.2,  p.87 and Theorem~3.7.3, p.96) and Rockafellar~\cite{MR43:445} (Corollary~10.1.1, p.83, Theorem~10.4,  p.86 and Theorem~25.3, p.244). Finally, we recall that $f$ is a concave function if $-f$ is a convex function.
\medskip
\\
Let $(\Omega,\Upsilon,\mu)$ be a measure space and let $\omega_1,\omega_2:\Omega\tends[0,\infty)$ be two $\mu$--measurable weights. In order to establish our generalized Hölder's inequality, we need the following hypotheses.
\begin{gather}
 \left\{
  \begin{split}
   \label{hypo1}
    & \Phi:I_1\tends[0,\infty) \mbox{ is a concave function, } I_1 \mbox{ is an} \\
    & \mbox{open interval and for a.e. } x\in\Omega, \; \omega_1(x)\in I_1,
  \end{split}
 \right.
\end{gather}
\begin{gather}
 \left\{
  \begin{split}
   \label{hypo2}
    &\Psi:I_2\tends[0,\infty) \mbox{ is a concave function, } I_2 \mbox{ is an} \\
    & \mbox{open interval and for a.e. } x\in\Omega, \; \omega_2(x)\in I_2,
  \end{split}
 \right.
\end{gather}
\begin{gather}
 \label{hypo3}
  \mbox{for a.e. } x\in\Omega,\; 1\le\Phi(\omega_1(x))\Psi(\omega_2(x)).
\end{gather}

\begin{thm}
\label{ineg1}
Let $(\Omega,\Upsilon,\mu)$ be a measure space, $\omega_1,\omega_2:\Omega\tends[0,\infty)$ be two  $\mu$--measurable weights and  $0<p<\infty.$ If there exist two functions $\Phi$ et $\Psi$ satisfying $(\ref{hypo1})-(\ref{hypo3})$ then for any $f\in L^p(\Omega,\Upsilon,\mu),$ $f\not\equiv0,$ we have
\begin{gather}
\label{ineg1-1}
1\le\Phi\left(\frac{\dsp\vint_\Omega|f|^p\omega_1d\mu}{\|f\|^p_{L^p(\Omega,\Upsilon,\mu)}}\right)
\Psi\left(\frac{\dsp\vint_\Omega|f|^p\omega_2d\mu}{\|f\|^p_{L^p(\Omega,\Upsilon,\mu)}}\right),
\end{gather}
as soon as $f\in L^p(\Omega,\Upsilon,d\mu)\cap L^p(\Omega,\Upsilon,\omega_1d\mu)\cap L^p(\Omega,\Upsilon,\omega_2d\mu).$
\end{thm}

\noindent
Obviously, one of the main issues to be clarified is whether there exist functions $\Phi$ and $\Psi$ satisfying the requirements (\ref{hypo1}), (\ref{hypo2}) and (\ref{hypo3}). This, of course, depends on the properties that the weight functions $\omega_1$ and $\omega_2$ satisfy. Below we shall give sufficient conditions on the weights $\omega_1$ and $\omega_2$ guaranteeing that $\Phi$ and $\Psi$ as above exist. This can be done by imposing some stronger conditions on the weight functions. More precisely, assume  that $\Omega=(m,\infty)$ (for some $m\in\R),$ $d\mu=dx$ is the Lebesgue's measure and
\begin{gather}
\label{hypo5}
\omega_1:(m,\infty)\tends(0,\omega_1(m)) \mbox{ is a convex and decreasing function and }
\omega_1(\infty)=0, \\
\label{hypo6}
\omega_2:(m,\infty)\tends(0,\infty) \mbox{ is a convex and increasing function and } \omega_2(\infty)=\infty, \\
\label{hypo7}
\Phi:(0,\omega_1(m))\tends(0,\infty) \mbox{ is a concave and increasing function and } \Phi(0)=0, \\
\label{hypo8}
\Psi:(\omega_2(m),\infty)\tends(0,\infty) \mbox{ is a concave and increasing function and } \Psi(\infty)=\infty, \\
\label{hypo9}
\forall t\in(m,\infty),\; 1\le\Phi(\omega_1(t))\Psi(\omega_2(t)).
\end{gather}

\noindent
Note that in (\ref{hypo7}), hypothesis $\Phi(0)=0$ means that $\Phi$ can be extended by continuity in $0$ by $0.$ 
\medskip
\\
The following result asserts that functions satisfying (\ref{hypo7})--(\ref{hypo9}) (and so (\ref{hypo1})--(\ref{hypo3})) exist, if the weights $\omega_1$ and $\omega_2$ verify the additional assumptions (\ref{hypo5})--(\ref{hypo6}).

\begin{thm}
\label{thmineg1}
Let $m>0$ and let $\omega_1,$ $\omega_2,$ be two weights satisfying $(\ref{hypo5})-(\ref{hypo6}).$ We define the function $\vphi$ by
\begin{gather}
\label{defpropineg1}
\forall t>m, \; \vphi(t)=m\frac{\omega_1(t)}{t}.
\end{gather}
Then the following assertions hold.
\begin{enumerate}
\item
\label{thmineg1-1}
The function $\Phi$ defined on $[0,\omega_1(m))$ by $\Phi(0)=0$ and $\Phi(t)=\dfrac{1}{\vphi^{-1}(t)},$ for  $t\not=0,$ satisfies $(\ref{hypo7}).$
\item
\label{thmineg1-2}
The function $\Psi$ defined on $(\omega_2(m),\infty)$ by $\Psi(t)=\omega_2^{-1}(t)$ satisfies $(\ref{hypo8}).$
\item
\label{thmineg1-3}
For $\Phi$ and $\Psi$ defined as above, estimate $(\ref{hypo9})$ holds.
\end{enumerate}
\end{thm}

\noindent
Before proving Theorems~\ref{ineg1}--\ref{thmineg1}, let us establish some preliminaries lemmas. The following result being a direct consequence of the definition of convex functions, we omit the proof.

\begin{lem}
\label{convexe}
Let $I\subset\R$ be an interval and let $\vphi:I\tends\R$ be a function. Then $\vphi$ is increasing and concave on $I$ if and only if $\vphi^{-1}$ is increasing and convex on $\vphi(I).$
\end{lem}

\noindent
The next lemma is the inverse version of the classical Jensen's inequality (W.~Rudin~\cite{MR924157}).

\begin{lem}[Inverse Jensen's inequality]
\label{jensen}
Let $(\Omega,\Upsilon,\nu)$ be a measure space such that $\nu(\Omega)=1$ and let $-\infty\le a<b\le+\infty.$ Assume that \\
$1)$ $\vphi:(a,b)\tends\R$ is a concave function, \\
$2)$ $f\in L^1(\Omega,\Upsilon,\nu)$ is such that for almost every $x\in\Omega,$ $f(x)\in(a,b).$ \\
Then $\vphi(f)_+\in L^1(\Omega,\Upsilon,\nu)$ and
\begin{gather}
\label{jensen1}
\vint_\Omega\vphi(f)d\nu\le\vphi\left(\:\vint_\Omega fd\nu\right).  
\end{gather}
\end{lem}

\begin{rmk}
\label{rmkjensen}
Since $\vphi$ is concave on $(a,b),$ it is continuous and $\vphi\circ f$ is a $\Upsilon$-measurable function. Furthermore, $\vphi(f)_+\in L^1(\Omega,\Upsilon,\nu)$ so the left-hand side of $(\ref{jensen1})$ makes sense and $\dsp\int_\Omega\vphi(f)d\nu\in[-\infty,+\infty).$ Indeed, since $\vphi$ is a concave function, it follows from the discussion at the beginning of this section that for any $t,s\in(a,b),$ $\vphi(t)\le\vphi(s)+\vphi^\prime_\ell(s)(t-s).$ In particular,
\begin{gather}
\label{rmkjensen1}
\vphi(f)\le\vphi(t_0)+\vphi^\prime_\ell(t_0)(f-t_0), \mbox{ a.e. in }\Omega, \medskip \\
\nonumber
\vphi(f)_+\le|\vphi(t_0)|+|\vphi^\prime_\ell(t_0)|(|f|+|t_0|)\in L^1(\Omega,\Upsilon,\nu),
\end{gather}
where $t_0=\dsp\int_\Omega fd\nu.$ Integrating \eqref{rmkjensen1} over $\Omega,$ we obtain \eqref{jensen1}. For more details, see Theorem~3.3 p.62 in W.~Rudin~\cite{MR924157}. 
\end{rmk}

\noindent
Now, we are in the conditions  to prove Theorem~\ref{ineg1}.
\medskip
\\
\begin{vproof}{of Theorem~\ref{ineg1}.}
Let $0<p<\infty,$ let $f\in L^p(\Omega,\Upsilon,\mu)\cap L^p(\Omega,\Upsilon,\omega_1d\mu)\cap L^p(\Omega,\Upsilon,\omega_2d\mu),$ $f\not\equiv0,$ and let $\nu$ be the measure defined by $\nu=\frac{|f|^p}{\|f\|^p_{L^p(\Omega,\Upsilon,\mu)}}\mu.$ Then $\nu(\Omega)=1.$ We apply twice Lemma~\ref{jensen} with $\vphi_1=\Phi,$ $f_1=\omega_1,$ $\vphi_2=\Psi$ and $f_2=\omega_2.$ Then $\Phi\circ\omega_1\in L^1(\Omega,\Upsilon,\nu),$ $\Psi\circ\omega_2\in L^1(\Omega,\Upsilon,\nu)$ and it follows from (\ref{hypo3}), Cauchy-Schwarz's inequality and (\ref{jensen1}) that
\begin{align*}
1=\left(\:\vint_\Omega1^\frac{1}{2}d\nu\right)^2
& \le\left(\:\vint_\Omega\Phi^\frac{1}{2}(\omega_1(x))\Psi^\frac{1}{2}(\omega_2(x))d\nu(x)\right)^2 \bigskip \\
& \le\vint_\Omega\Phi(\omega_1(x))d\nu(x)\vint_\Omega\Psi(\omega_2(x))d\nu(x) \bigskip \\
& \le\Phi\left(\:\dsp\vint_\Omega\omega_1(x)d\nu(x)\right)
\Psi\left(\:\dsp\vint_\Omega\omega_2(x)d\nu(x)\right) \bigskip \\
& =\Phi\left(\frac{\dsp\vint_\Omega|f|^p\omega_1d\mu}{\|f\|^p_{L^p(\Omega,\Upsilon,\mu)}}\right)
\Psi\left(\frac{\dsp\vint_\Omega|f|^p\omega_2d\mu}{\|f\|^p_{L^p(\Omega,\Upsilon,\mu)}}\right).
\end{align*}
Hence (\ref{ineg1-1}).
\medskip
\end{vproof}

\noindent
The proof of Theorem \ref{thmineg1} relies on the following lemma.

\begin{lem}
\label{concave}
Let $m\in[0,\infty),$ $0<M\le\infty$ and $p\in[1,\infty).$ Let $f:(m,\infty)\tends(0,M)$ be a nonincreasing function such that $f(m)=M.$ Define the function $\vphi_p$ on $(m,\infty)$ by
\begin{gather}
\label{concave1}
\forall t>m, \; \vphi_p(t)=\frac{f(t)}{t^p}.
\end{gather}
If $f$ is convex on $(m,\infty)$ then $\vphi_p$ is convex on $(m,\infty)$ and $\dfrac{1}{\vphi_p^{-1}}$ is  concave and increasing on $\left(0,\frac{M}{m^p}\right),$ where we have used the notation $\frac{M}{m^p}=+\infty$ if $m=0$ and/or $M=+\infty.$ Furthermore, $\vlim_{t\searrow0}\frac{1}{\vphi^{-1}(t)}=0.$ 
\end{lem}

\begin{rmk}
If $0<p<1$ then the conclusion of Lemma~\ref{concave} may be false. Indeed, let $q_0\in(p,1)$ and set $q=\frac{1}{q_0}>1.$ We then choose $f(t)=\frac{1}{t^{q_0-p}},$ $t>0.$ Then $f$ and $\vphi_p$ are obviously convex and decreasing on $(0,\infty).$ But for any $t>0,$ $\frac{1}{\vphi_p^{-1}(t)}=t^q.$ So that $\vphi_p$ is not concave on $(0,\infty)$ since $q>1.$
\end{rmk}

\begin{rmk}
\label{rmkpropphi}
Let $f:(m,\infty)\tends(0,\infty)$ be an application, where $m\in\R.$ Assume that $f$ is convex on $(m,\infty)$ and that $\vlim_{t\to\infty}f(t)=0.$ If $f$ is nonincreasing on $(m,\infty)$ then it is in fact decreasing on $(m,\infty).$ Indeed, if $f$ is not decreasing on $(m,\infty)$ then $f(t)=f(a)>0$ for any $t\in(a,b),$ for some interval $(a,b)\subset(m,\infty).$ Since $\vlim_{t\to\infty}f(t)=0,$ we necessarily have $b<\infty.$ Then $f^\prime\equiv0$ on $(a,b)$ and, by hypothesis $\vlim_{t\to\infty}f(t)=0,$ this implies that $f^\prime(t_0)<0,$ for some $t_0\in(b,\infty).$ This contradicts hypothesis $f$ is convex.
\end{rmk}

\begin{vproof}{of Lemma \ref{concave}.}
Let $\vphi_p$ be defined by (\ref{concave1}). Note that $\vphi_p:\left(m,\infty\right)\tends\left(0,\frac{M}{m^p}\right)$ being bijective, continuous and decreasing, it follows that
$\vphi_p^{-1}:\left(0,\frac{M}{m^p}\right)\tends\left(m,\infty\right)$ is well-defined, continuous and decreasing. So $\frac{1}{\vphi_p^{-1}}:\left(0,\frac{M}{m^p}\right)\tends\left(0,\frac{1}{m}\right)$ is continuous and increasing, where we have used the notation $\frac{1}{m}=+\infty$ if $m=0.$ The product of two positive and convex functions with the same monotonicity being convex, it follows that the function $t\longmapsto\frac{f(t)}{t^p}$ is convex and so $\vphi_p$ is convex. Moreover, hypothesis $\vlim_{t\nearrow\infty}\vphi(t)=0$ implies that $\vlim_{t\searrow0}\frac{1}{\vphi^{-1}(t)}=0.$ Since $f$ is convex, according to the basic properties on convex functions we recalled in the beginning of this section, there exists a sequence $(a_n)_{n\in\N}\subset(m,\infty)$ such that $f$ is $\vC^1$ and $f^\prime$ is nondecreasing relatively to $(m,\infty)\setminus\vN,$ with $\vN=\bigcup\limits_{n=1}^\infty\{a_n\}.$ Now, we proceed to the proof in 3 steps. \\
{\bf Step 1.} Set for every $t\in(m,\infty)\setminus\vN,$
\begin{gather}
\label{demconcave1}
h(t)=-(f^\prime(t)t-pf(t)) \quad \mbox{ and } \quad
g(t)=\dfrac{h(t)}{t^{p-1}}.
\end{gather}
Then $g$ is nonincreasing and nonnegative on $(m,\infty)\setminus\vN.$ \\
Indeed, let $s,t\in(m,\infty)\setminus\vN$ be such that $s<t.$ Since $f$ is convex, it follows from the discussion at the beginning of this section that $f(t)-f(s)\le f^\prime(t)(t-s).$ Using this estimate, $p\ge1$ and again the fact that $f$ is nonincreasing and $f^\prime$ is nondecreasing relatively to $(m,\infty)\setminus\vN,$ we obtain that
\begin{align*}
        & h(t)-h(s)=p(f(t)-f(s))-(t-s)f^\prime(t)-s(f^\prime(t)-f^\prime(s)) \medskip \\
\le \; & f(t)-f(s)-f^\prime(t)(t-s)\le0.
\end{align*}
Consequently, $h$ is is nonincreasing. Since it is nonnegative (because $f$ is nonnegative and nonincreasing), it follows that $g$ is also nonincreasing and nonnegative relatively to $(m,\infty)\setminus\vN.$ \\
{\bf Step 2.} We claim that, for any $t>m,$ $\vphi_p(t)=\dsp\int_0^{1/t}g\left(\frac{1}{s}\right)ds.$ \\
Indeed, by (\ref{concave1})--(\ref{demconcave1}), we have for every $\sigma\in(m,\infty)\setminus\vN,$
\begin{gather*}
-\vphi_p^\prime(\sigma)=-\dfrac{f^\prime(\sigma)\sigma^p-pf(\sigma)\sigma^{p-1}}{\sigma^{2p}}
=-\dfrac{f^\prime(\sigma)\sigma-pf(\sigma)}{\sigma^{p+1}}
=\dfrac{h(\sigma)}{\sigma^{p-1}}\dfrac{1}{\sigma^2}=\dfrac{g(\sigma)}{\sigma^2}.
\end{gather*}
Then for any $\eps>0,$ $\vphi_p^\prime\in L^1(m+\eps,\infty)$ and so $\vphi_p(t)=\dsp\int_t^\infty\frac{g(\sigma)}{\sigma^2}d\sigma,$ which yields the desired result, by using the change of variables $\sigma=\frac{1}{s}.$ \\
{\bf Step 3.} Conclusion. \\
Let $\psi$ be defined on $\left(0,\frac{M}{m^p}\right)$ by $\psi(t)=\frac{1}{\vphi_p^{-1}(t)}.$ Thus by Step~2, we have for any $t\in\left(0,\frac{1}{m}\right),$
\begin{gather*}
 \psi^{-1}(t)=\vphi_p\left(\frac{1}{t}\right)=\vint_0^t g\left(\frac{1}{s}\right)ds.
\end{gather*}
Then $\psi^{-1}$ is absolutely continuous and for almost every $t\in\left(0,\frac{1}{m}\right),$
$\left(\psi^{-1}\right)^\prime(t)=g\left(\frac{1}{t}\right)\ge0.$ Since $g$ is nonincreasing relatively to $(m,\infty)\setminus\vN$ (Step~1), it follows that $\psi^{-1}$ is increasing and convex on $\left(0,\frac{1}{m}\right).$ By Lemma~\ref{convexe}, $\psi\stackrel{{\rm def}}{=}\frac{1}{\vphi_p^{-1}}$ is increasing and concave on $\left(0,\frac{M}{m^p}\right).$ Hence the result.
\medskip
\end{vproof}

\begin{vproof}{of Theorem \ref{thmineg1}.}
Let $\vphi$ be defined on $(m,\infty)$ by (\ref{defpropineg1}). By (\ref{hypo5})--(\ref{hypo6}), $\omega_2$ is invertible on $(\omega_2(m),\infty)$ and $\vphi:(m,\infty)\tends(0,\omega_1(m))$ is a bijective and decreasing function. Then definition of $\Phi$ and $\Psi$ makes sense. \\
{\bf Proof of \ref{thmineg1-1}--\ref{thmineg1-2}.}
Assertion~\ref{thmineg1-1} is a direct consequence of Lemma~\ref{concave} applied to $f=m\omega_1$ and assertion~\ref{thmineg1-2} comes from \eqref{hypo6} and Lemma~\ref{convexe}. \\
{\bf Proof of \ref{thmineg1-3}.}
By (\ref{defpropineg1}) and definition of $\Phi,$ $\Phi^{-1}\left(\frac{1}{t}\right)=\vphi(t)\le\omega_1(t),$ for any $t>m.$ Since $\vphi$ and $\omega_1$ are both decreasing, this implies that
$$
\forall t\in(0,\omega_1(m)), \; \Phi(t)=\frac{1}{\vphi^{-1}(t)}\ge\frac{1}{\omega_1^{-1}(t)}.
$$
With the above estimate, we obtain that
$$
\forall t>m, \;
\Phi(\omega_1(t))\Psi(\omega_2(t))=\Phi(\omega_1(t))t\ge\frac{t}{\omega_1^{-1}(\omega_1(t))}=1.
$$
Hence (\ref{hypo9}). This concludes the proof.
\medskip
\end{vproof}

\noindent
We now give an example where the assumptions of Theorem~\ref{ineg1} are satisfied. The weight functions $\omega_1$, $\omega_2$ are of a particular form that arises naturally in applications: While $\omega_1$ tends to zero exponentially at $\infty$, $\omega_2$ grows as a polynomial function. This is a case that may not be covered by H\"older's inequality. In the sequel, we compute explicitly the functions $\Phi$ and $\Psi$ for which the generalized interpolation inequality holds.

\begin{exa}
\label{ineg2}
Let $\Omega=\R^N\setminus \ovl B(0,1)$ and $A\ge1.$ We consider the weights defined on $\Omega$ by $\omega_1(x)=e^{-A|x|}$ and $\omega_2(x)=|x|^2.$ We define the interpolating functions $\Psi(t)=\sqrt{t}$ $(t\ge0)$ and
$$
\begin{array}{rl}
 \forall t\in[0,e^{A-2}], & \Phi(t)=\left\{
  \begin{array}{rl}
                                0, & \mbox{if } \; t=0, \medskip \\
   \dfrac{2A}{A-\ln t}, & \mbox{if } \; 0<t\le e^{A-2}.
  \end{array}
 \right.
\end{array}
$$
The hypotheses of Theorem 2.1 are satisfied since the weights $\omega_1$ and $\omega_2$ and the interpolation functions $\Phi$ and $\Psi$ defined as above, satisfy the pointwise inequality (\ref{hypo3}) as it is immediate to check. Indeed, for any $x\in\Omega,$
$$
\Phi(\omega_1(x))\Psi(\omega_2(x))=\frac{2A|x|}{A+A|x|}=\frac{2|x|}{1+|x|}\ge1,
$$
since $|x|>1.$ Moreover, a straightforward calculation shows that $\Phi$ is concave on $[0,e^{A-2}].$ As a consequence of Theorem 2.1 we obtain the following functional generalized interpolation inequality. Let $f\in L^2(\Omega;\C)\setminus\{0\}$ be such that $|\: . \:|f(\: . \:)\in L^2(\Omega;\C).$ Then,
\begin{gather}
\label{corineg1}
\|f\|_{L^2(\Omega)}\le2\sqrt{\vint_\Omega|f(x)|^2|x|^2dx}
\dfrac{A}{A+\ln\left(\dfrac{1}{\|f\|_{L^2(\Omega)}^2}\dsp\vint_\Omega|f(x)|^2e^{-A|x|}dx\right)}.
\end{gather}
In the same way, we have
\begin{gather}
\label{corineg2}
\|u\|_{\ell^2(\N)}\le2\sqrt{\dsp\vsum_{n=1}^\infty n^2|u_n|^2}
\dfrac{A}{A-\ln\left(\dfrac{1}{\|u\|_{\ell^2(\N)}^2}\dsp\vsum_{n=1}^\infty e^{-An}|u_n|^2\right)},
\end{gather}
for any $u=(u_n)_{n\in\N}\in\ell^2(\N;\C)\setminus\{0\}$ such that $(nu_n)_{n\in\N}\in\ell^2(\N;\C).$ Note that one always has for any $A\ge1,$
$$
0<\frac{1}{\|f\|_{L^2(\Omega)}^2}\dsp\vint_\Omega|f(x)|^2e^{-A|x|}dx\le e^{-A}\le e^{A-2}
$$
and
$$
0<\frac{1}{\|u\|_{\ell^2(\N)}^2}\dsp\vsum_{n=1}^\infty e^{-An}|u_n|^2\le e^{-A}\le e^{A-2},
$$
(since $e^{-A}\le e^{A-2}\iff A\ge1)$ so the above quantities takes their values in the domain of concavity of $\Phi.$ It follows that estimates (\ref{corineg1}) and (\ref{corineg2}) always make sense.
\end{exa}

\section{Optimality}
\label{optimality}

It this section, we discuss the notion of optimality for the pairs of functions $(\Phi,\Psi)$ satisfying the interpolation inequalities above. We will also give sufficient conditions guaranteeing the pair is optimal. Throughout this section,  for simplicity, we assume that $\Omega=(m,\infty)$ (for some $m\in\R)$ and that $d\mu=dx$ is the Lebesgue's measure. Before introducing the definition of optimality, we need the following lemma.

\begin{lem}
\label{prelem}
Let $m\in\R$ and let $\omega_1,$ $\omega_2,$ $\Phi$ and $\Psi$ satisfy $(\ref{hypo5})-(\ref{hypo9}).$ Let $\delta\in(0,\omega_1(m)]$ be such that $\Phi(\delta)= \frac{1}{\Psi\left(\omega_2(m)\right)},$ if $\Psi\left(\omega_2(m)\right)>0$ and let $\delta=+\infty,$ if $\Psi\left(\omega_2(m)\right)=0.$ We define
\begin{gather}
\label{optim1}
\forall t\in(0,\delta), \; \vH_{\Phi,\Psi}(t)=\frac{1}{\Psi^{-1}\left(\dfrac{1}{\Phi(t)}\right)}.
\end{gather}
Then $\vH_{\Phi,\Psi}$ is a positive, increasing and continuous function on $(0,\delta)$ and $\vlim_{t\searrow0}\vH_{\Phi,\Psi}(t)=0.$ Furthermore,
\begin{gather}
\label{prelem-1}
\forall t\in(0,\delta), \;
0<\dfrac{1}{\omega_2\circ\omega_1^{-1}(t)}\le\vH_{\Phi,\Psi}(t).
\end{gather}
Finally,
\begin{gather}
\label{prelem-2}
\vH_{\Phi,\Psi}^{-1}(t)=\Phi^{-1}\left(\frac{1}{\Psi\left(\frac{1}{t}\right)}\right),
\end{gather}
for any $t\in\left(0,\vH_{\Phi,\Psi}(\delta)\right).$
\end{lem}

\begin{rmk}
\label{rmkthmoptim-1}
Note that such a $\delta\in (0, \omega_1(m)]$   exists because of the continuity of $\Phi$.
\end{rmk}

\noindent
Assuming for the moment that Lemma~\ref{prelem} holds (we shall return to its proof later), the following definition makes sense.

\begin{defi}
\label{defioptim}
Let $m\in\R$ and $\omega_1,$ $\omega_2,$ $\Phi$ and $\Psi$ satisfy (\ref{hypo5})--(\ref{hypo9}). We say that ($\Phi$, $\Psi$)  is an {\it optimal pair  for the weights} $(\omega_1, \omega_2)$  if the function $\vH_{\Phi,\Psi}$ defined by (\ref{optim1}) satisfies
\begin{gather}
\label{defoptim1-1}
\vH_{\Phi,\Psi}\stackrel{0}{\thickapprox}\dfrac{1}{\omega_2\circ\omega_1^{-1}}.
\end{gather}
Here and in the sequel, by $\vH_{\Phi,\Psi}\stackrel{0}{\thickapprox}\dfrac{1}{\omega_2\circ\omega_1^{-1}}$ we mean that there exist two constants $C>0$ and $\eps\in(0,\delta)$ such that
\begin{gather}
\label{equivfunc-1}
\forall t\in(0,\eps), \;
\dfrac{1}{\omega_2\circ\omega_1^{-1}(t)}\le\vH_{\Phi,\Psi}(t)\le\dfrac{C}{\omega_2\circ\omega_1^{-1}(t)},
\end{gather}
where $\delta>0$ is given in Lemma~\ref{prelem}.
\end{defi}

\noindent
In view of (\ref{prelem-1}) when (\ref{defoptim1-1}) holds, the function $\vH_{\Phi,\Psi}(t)$ goes to $0$ as $t\searrow0$ as rapidly as possible. The pair $(\Phi, \Psi)$ is then optimal in that sense. As we shall see in applications, this will yield the optimal decay rate for the energy of solutions of damped wave-like equations.

\begin{rmk}
\label{rmkdef}
It is important to note that the notion of optimal pair $(\Phi, \Psi)$ depends on the  weights $(\omega_1, \omega_2)$. On the other hand, given two weights  $\omega_1$ and $\omega_2$  satisfying (\ref{hypo5})--(\ref{hypo6}) and  a pair $(\Phi,\Psi)$ satisfying (\ref{hypo8})--(\ref{hypo9}), if $\Phi^{-1}\left(\frac{1}{\Psi\circ\omega_2}\right)$ is convex then the pair $(\Phi,\Psi)$ is necessarily optimal with respect  to the weights $\widetilde{\omega_1}$ and $\omega_2,$ where we have chosen $\widetilde{\omega_1}(t)=\Phi^{-1}\left(\frac{1}{\Psi(\omega_2(t))}\right).$ Indeed, (\ref{hypo5})--(\ref{hypo8}) hold for $(\widetilde{\omega_1},\omega_2,\Phi,\Psi).$ Moreover, 
$$
\Phi(\widetilde{\omega_1}(t))\Psi(\omega_2(t))=\frac{1}{\Psi(\omega_2(t))}\Psi(\omega_2(t))=1,
$$
and (\ref{hypo9}) is fulfilled. Finally, a straightforward calculation gives
$$
\vH_{\Phi,\Psi}(t)\stackrel{\mathrm{def}}{=}\frac{1}{\Psi^{-1}\left(\dfrac{1}{\Phi(t)}\right)}
=\dfrac{1}{\omega_2\circ\widetilde{\omega_1}^{-1}(t)}.
$$
Hence (\ref{defoptim1-1}).
\end{rmk}

\noindent
Now we give a sufficient condition for the pair $(\Phi, \Psi)$ to be optimal.

\begin{prop}
\label{propoptim1}
Let $m\in\R$ and let $\omega_1$ and $\omega_2$ be satisfying $(\ref{hypo5})-(\ref{hypo6}).$ Let $1\le p<\infty,$ and set  
\begin{gather}
\label{psi}
\forall t>\omega_2(m), \quad\Psi_p(t)=\left(\omega_2^{-1}(t)\right)^\frac{1}{p}, 
\end{gather}
and 
\begin{gather}
\label{phi}
\forall t\in(0,\omega_1(m)), \quad \Phi_p(t)=\dfrac{1}{\left(\omega_1^{-1}(t)\right)^\frac{1}{p}},
\end{gather}
together with $\Phi_p(0)=0$.

If $\dfrac{1}{\left(\omega_1^{-1}\right)^\frac{1}{p}}$ is concave on $(0,\omega_1(m))$ then $(\Phi_p,\Psi_p)$ constitutes an optimal pair for the weights $(\omega_1, \omega_2)$.
\end{prop}

\noindent
On the other hand, the following Proposition guarantees that, once we have an optimal pair $(\Phi, \Psi)$ it is easy to build other optimal pairs. Of course, in practice, when applying the interpolation inequalities to obtain decay rates for evolution equations, it is irrelevant whether one uses an optimal pair or another since all of them, by definition, yield the same decay rates.

\begin{prop}
\label{propoptim2}
Let $m\in\R$ and let $\omega_1,$ $\omega_2,$ $\Phi$ and $\Psi$ be satisfying $(\ref{hypo5})-(\ref{hypo7}).$ Let $0<p<\infty,$ let $(0,\delta)$ be the interval of definition of $\vH_{\Phi,\Psi}$ and let $(0,\delta_p)$ be the interval of definition of $\vH_{\Phi^p,\Psi^p}$ $($see Lemma~$\ref{prelem}).$ Then
\begin{gather*}
\forall t\in(0,\inf\{\delta,\delta_p\}), \; \vH_{\Phi,\Psi}(t)=\vH_{\Phi^p,\Psi^p}(t).
\end{gather*}
In particular, if $(\Phi, \Psi)$ is an optimal pair for the weights $(\omega_1, \omega_2)$, then the same holds for $(\Phi^p,\Psi^p).$
\end{prop}

\begin{rmk}
In other words, Proposition~\ref{propoptim2} means that, from the point of view of the decay of $\vH_{\Phi,\Psi},$ the inequalities
$
1\le\Phi(\omega_1)\Psi(\omega_2)$ and $
1\le\Phi^p(\omega_1)\Psi^p(\omega_2),
$ yield the same result.
\end{rmk}

\begin{vproof}{of Lemma \ref{prelem}.}
Let $\Phi$ and $\Psi$ be any functions satisfying (\ref{hypo8})--(\ref{hypo9}) and $\delta>0$ be defined as in Lemma~\ref{prelem}. It follows from (\ref{hypo5})--(\ref{hypo9}) and definition of $\delta$ that
$$
\forall t\in(0,\omega_1(m)), \;
1\le\Phi(t)\Psi\left(\omega_2\circ\omega_1^{-1}(t)\right)
\quad \mbox{ and } \quad
\forall t\in(0,\delta), \;
0<\Phi(t)<\frac{1}{\Psi(\omega_2(m))}\le+\infty.
$$
We then have
$$
\forall t\in(0,\delta), \;
0\le\Psi(\omega_2(m))<\dfrac{1}{\Phi(t)}\le\Psi\left(\omega_2\circ\omega_1^{-1}(t)\right).
$$
Since $\Psi^{-1}$ is increasing on $(\Psi(\omega_2(m)),\infty),$ this gives
$$
\forall t\in(0,\delta), \;
0<\Psi^{-1}\left(\dfrac{1}{\Phi(t)}\right)\stackrel{{\rm def}}{=}
\dfrac{1}{\vH_{\Phi,\Psi}(t)}\le\omega_2\circ\omega_1^{-1}(t),
$$
which yields (\ref{prelem-1}). Properties of $\vH_{\Phi,\Psi}$ follows easily from (\ref{hypo7})--(\ref{hypo8}).
\medskip
\end{vproof}

\begin{vproof}{of Proposition~\ref{propoptim2}.}
Let $s\in\vH_{\Phi,\Psi}((0,\delta))\cap\vH_{\Phi^p,\Psi^p}((0,\delta_p)).$ Then we have,
$$
\begin{array}{rccccl}
& \vH_{\Phi^p,\Psi^p}(t)=s & \iff & \dfrac{1}{\left(\Psi^p\right)^{-1}\left(\dfrac{1}{\Phi^p(t)}\right)}=s & \iff &
\left(\Psi^p\right)^{-1}\left(\dfrac{1}{\Phi^p(t)}\right)=\dfrac{1}{s} \medskip\\
\iff & \dfrac{1}{\Phi^p(t)}=\Psi^p\left(\dfrac{1}{s}\right) & \iff & \dfrac{1}{\Phi(t)}=\Psi\left(\dfrac{1}{s}\right) & \iff &
\vH_{\Phi,\Psi}(t)=s.
\end{array}
$$
Hence the result.
\medskip
\end{vproof}

\begin{vproof}{of Proposition \ref{propoptim1}.}
Assume that hypotheses of Proposition~\ref{propoptim1} are satisfied. It follows from Lemma~\ref{convexe} and (\ref{hypo6}) that $\Psi_p$ satisfies (\ref{hypo8}). By (\ref{hypo5}) and the fact that $\frac{1}{\left(\omega_1^{-1}\right)^\frac{1}{p}}$ is concave on $(0,\omega_1(m)),$ the function $\Phi_p$ defined as in (\ref{phi}) satisfies (\ref{hypo7}). By (\ref{psi}) and (\ref{phi}), (\ref{hypo9}) and (\ref{defoptim1-1}) are verified. 
Indeed, by Proposition~\ref{propoptim2},
\begin{align*}
& \vH_{\Phi_p,\Psi_p}(t)=\vH_{\Phi_p^p,\Psi_p^p}(t)=\frac{1}{(\Psi_p^p)^{-1}\left(\dfrac{1}{\Phi_p^p(t)}\right)} \\
= & \frac{1}{(\Psi_p^p)^{-1}(\omega_1^{-1}(t))}=\frac{1}{(\omega_2^{-1})^{-1}(\omega_1^{-1}(t))}
=\frac{1}{\omega_2\circ\omega_1^{-1}(t)}.
\end{align*}
This concludes the proof.
\medskip
\end{vproof}

\begin{rmk}
\label{rmkthmoptim0}
Note that the hypothesis $p\ge1$ in Proposition~\ref{propoptim1} is made to ensure that $\left(\omega_2^{-1}\right)^\frac{1}{p}$ is a concave function. So it follows from the above proof that, if $0<p<1$ is such that  $\left(\omega_2^{-1}\right)^\frac{1}{p}$ is  concave, then the conclusion of Proposition~\ref{propoptim1} still holds.
\end{rmk}

\begin{rmk}
\label{rmkthmoptim3}
Proposition~\ref{propoptim2} shows the non uniqueness of the optimal pairs $(\Phi, \Psi)$. One may give other examples. Let $m\in\R$ and let $\omega_1$ and $\omega_2$ be satisfying (\ref{hypo5})--(\ref{hypo6}). Following the proof of Proposition~\ref{propoptim1}, we can show that if $\frac{1}{\omega_2\circ\omega_1^{-1}}$ is concave then the functions $\Psi=\Id$ and $\Phi=\frac{1}{\omega_2\circ\omega_1^{-1}}$ are an optimal pair of functions.
\end{rmk}

\section{Application to the stabilization on the wave equation with Dirichlet boundary condition}
\label{damping}

In this section, we give some applications of Section~\ref{interpolation}. We recover and extend the results of Ammari, Henrot and Tucsnak~\cite{MR2002j:93073}, Ammari and Tucsnak~\cite{MR1814271} and Jaffard, Tucsnak and Zuazua~\cite{MR99g:93073}. We will detail the first example (Subsection~\ref{sswave}) and we will indicate how we proceed for the others equations (for conciseness of the paper, we will not detail the proof, the method being very technical). We apply our interpolation inequality to the stabilization of a wave equation with a damping control concentrated on an interior point (Subsection~\ref{sswave}) and to the stabilization of a Bernoulli--Euler beam with a damping control concentrated in an interior point (Subsection~\ref{ssbeam}).

\subsection{Explanation of the method}
\label{explanation}

\numberwithin{equation}{subsection}
\renewcommand{\theequation}{\thesubsection.\arabic{equation}}

To set the context, we introduce some notations and refer to Ammari and Tucsnak~\cite{MR2002f:93104} for more details. We consider $u$ the solution of the following equation.
\begin{gather}
 \left\{
  \begin{split}
   \label{eqrefu}
    u_{tt}+Au+BB^\star u_t=0, & \quad (t,x)\in(0,\infty)\times I, \\
                  u(0,x)=u^0(x), & \quad x\in I, \\
               u_t(0,x)=u^1(x), & \quad x\in I,
   \end{split}
 \right.
\end{gather}
where $A$ is a linear unbounded self-adjoint operator, $B\in{\cal L}(U;D(A^\frac{1}{2})^\star),$ $(U,\|\:.\:\|_U)$ is a complex Hilbert space, $D(A^\frac{1}{2})=\ovl{D(A)}^{\|\: . \:\|_{\frac{1}{2}}},$ $\|u\|_{\frac{1}{2}}=\sqrt{\langle Au,u\rangle},$ $D(A^\frac{1}{2})^\star$ is the topological dual of the space $D(A^\frac{1}{2}),$ $I=(0,L)$ is an interval of $\R$ and where the initial data $(u^0,u^1)$ are chosen in a Banach space $V\times L^2(I),$ in which equation (\ref{eqrefu}) is well set. The associated energy $E$ of $u$ is given by
\begin{gather}
\label{E}
\forall t\ge0, \;
E(u(t))=\frac{1}{2}\left(\|u_t(t)\|_{L^2(I)}^2+\|A^\frac{1}{2}u(t)\|_{L^2(I)}^2\right),
\end{gather}
and satisfies
\begin{gather}
\label{decrefE}
\forall t\ge s\ge0, \;
E(u(t))-E(u(s))=-\vint_s^{\;t}\|(B^\star u)_t(\sigma)\|_U^2d\sigma\le0.
\end{gather}
Typically, $V\times L^2(I)=D(A^\frac{1}{2})\times L^2(I)$ is the space for which the energy is well-defined and $U=\R.$ But we need more regularity and we choose $(u^0,u^1)\in\vD(\vA),$ where
$$
\begin{array}{rl}
\vA= & \!\!\!\! \left(
 \begin{array}{cc}
    0 & \rm{Id}    \\
  -A & -BB^\star
 \end{array}
\right).
\end{array}
$$
Denote by $(a_n)_{n\ge0}$ the sequence of the Fourier's coefficient of $u^0$ and by $(b_n)_{n\ge0}$ the  $u^1$ one. We also consider $v$ the solution of
\begin{gather}
 \left\{
  \begin{split}
   \label{eqrefv}
          v_{tt}+Av=0, & \quad (t,x)\in(0,\infty)\times I, \\
       v(0,x)=u^0(x), & \quad x\in I, \\
    v_t(0,x)=u^1(x), & \quad x\in I.
   \end{split}
 \right.
\end{gather}
Depending of the spaces $V$ and $\vD(\vA)$ we have chosen, we obtain for $(u^0,u^1)\in D(A)\times V,$
\begin{gather*}
\|(u^0,u^1)\|_{\vD(\vA)}^2=\dsp\vsum_{n=0}^\infty n^p(a_n^2+b_n^2)\omega_2(n), \medskip \\
E(u(0))\stackrel{{\rm def}}{=}\frac{1}{2}\|(u^0,u^1)\|_{V\times L^2(I)}^2
=\frac{1}{2}\dsp\vsum_{n=0}^\infty n^p(a_n^2+b_n^2),
\end{gather*}
for some weight $\omega_2$ satisfying (\ref{hypo6}) and some $p\in[0,\infty).$ Roughly speaking, in our examples, this comes from the expansion of $u^0$ and $u^1$ in Fourier's series and Parseval's identity.
\medskip
\\
First, we show that there exist a time $T>0,$ two constants $C>0$ and $C_1>0$ and a weight $\omega_1$ satisfying (\ref{hypo5}), such that for any initial data $(u^0,u^1)\in V\times L^2(I),$
\begin{gather}
\label{refestiminf}
\vint_0^{\;T}\|(B^\star u)_t(t)\|_U^2dt\ge C\vint_0^{\;T}\|(B^\star v)_t(t)\|_U^2dt
\ge C_1\vsum_{n=0}^\infty n^p(a_n^2+b_n^2)\omega_1(n),
\end{gather}
where the last estimate comes from  Ingham's inequality (Ingham~\cite{MR1545625}). For a complete example, see Lemmas~\ref{equi} and \ref{borneinfv}.
\medskip
\\
Second, we define the weak energy $E_-$ and the strong energy $E_+$ as follow.
\begin{gather}
\label{refE+}
E_+(0)=\dsp\vsum_{n=0}^\infty n^p(a_n^2+b_n^2)\omega_2(n), \medskip \\
\label{refE}
E(0)=\dsp\vsum_{n=0}^\infty n^p(a_n^2+b_n^2), \medskip \\
\label{refE-}
E_-(0)=\dsp\vsum_{n=0}^\infty n^p(a_n^2+b_n^2)\omega_1(n).
\medskip
\end{gather}
Third, we show that there exist two functions $\Phi$ and $\Psi$ satisfying (\ref{hypo7}) and (\ref{hypo8}). From Theorem~\ref{ineg1}, we have (\ref{ineg1-1}). Typically, we choose $\Phi (t)=\frac{1}{\vphi^{-1}(t)}$ and $\Psi(t)=\omega_2^{-1}(t),$ where $\vphi(t)=\frac{\omega_1(t)}{t^p}$ with $p\in\{0,2,4\}.$ From (\ref{ineg1-1}) and (\ref{refE+})--(\ref{refE-}), we deduce that
\begin{gather}
\label{refestiminf2}
E_-(0)\ge E(0) \Phi^{-1} \left(\frac{1}{\Psi\left(\frac{E_+(0)}{E(0)}\right)}\right)
=E(0)\vH^{-1}_{\Phi,\Psi}\left(\dfrac{E(0)}{E_+(0)}\right),
\end{gather}
where $\vH_{\Phi,\Psi}^{-1}$ is defined by (\ref{prelem-2}). Putting together (\ref{decrefE}), (\ref{refestiminf}) and (\ref{refestiminf2}), we obtain
\begin{gather}
\label{refestiminf3}
E(T)\le E(0)-C_1E(0)\vH^{-1}_{\Phi,\Psi}\left(\dfrac{E(0)}{E_+(0)}\right).
\end{gather}
See Lemma~\ref{estiminf} for a complete example.
\medskip
\\
Fourth, we use (\ref{refestiminf3}), the semigroup property and the method of Ammari and Tucsnak~\cite{MR2002f:93104} to show that
\begin{gather}
  \forall t\ge0, \; E(t)\le
  C\vH_{\Phi,\Psi}\left(\dfrac{1}{t+1}\right)\|(u^0,u^1)\|_{\vD(\vA)}^2.
\end{gather}
Their proof is based on an interpolation method. See Theorem~\ref{thmdecay1} for a complete example.

\subsection{Notations for the wave equation (\ref{wave}) with Dirichlet boundary condition and known results}
\label{sswave}

We consider a wave equation with a damping control concentrated on an interior point $a\in(0,1)$ with homogenous Dirichlet boundary condition,
\begin{gather}
 \left\{
  \begin{split}
   \label{wave}
    u_{tt}-u_{xx}+\delta_{a}u_t(t,a)=0, & \quad (t,x)\in(0,\infty)\times(0,1), \\
       u(0,x)=u^0(x), \; u_t(0,x)=u^1(x), & \quad x\in(0,1), \\
                                       u(t,0)=u(t,1)=0 & \quad t\in[0,\infty).
   \end{split}
 \right.
\end{gather}

\noindent
Let $V_1=H^1_0(0,1).$ A direct calculation gives that for any $u\in V_1,$ $\|u\|_{L^2(0,1)}\le\|u_x\|_{L^2(0,1)},$ so we may endow $V_1$ of the norm $\|u\|_{V_1}=\|u_x\|_{L^2(0,1)},$ for any $u\in V_1.$ Let $X_1=V_1\times L^2(0,1),$
$$
Y_1=\left(H^1_0(0,1)\cap H^2(0,a)\cap H^2(a,1)\right)\times H^1_0(0,1),
\quad
D(A_1)=H^1_0(0,1)\cap H^2(0,1),
\quad
A_1=-\frac{{\mathrm d}^2}{{\mathrm d}x^2},
$$
\begin{gather*}
\vD(\vA_1)=\left\{(u,v)\in Y_1;\; \frac{{\mathrm d}u}{{\mathrm d}x}(a_+)
-\frac{{\mathrm d}u}{{\mathrm d}x}(a_-)=v(a)\right\},
\end{gather*}
with
$$
\|(u,v)\|_{\vD(\vA_1)}^2=\|(u,v)\|_{Y_1}^2=\|u\|_{H^2(0,a)}^2+\|u\|_{H^2(a,1)}^2+\|v\|_{H^1_0(0,1)}^2,
$$
and let
$
\begin{array}{rl}
\vA_1= & \!\!\!\! \left(
 \begin{array}{cc}
      0    &  \rm{Id}    \\
  -A_1  &  -\delta_{a}
 \end{array}
\right).
\end{array}
$
We define the energy $E_1$ for $u$ solution of equation (\ref{wave}) by
\begin{gather}
\label{energy1}
\forall t\ge0, \; E_1(u(t))=\frac{1}{2}\left(\|u_t(t)\|_{L^2(0,1)}^2 + \|u_x(t)\|_{L^2(0,1)}^2\right)
=\frac{1}{2}\|(u(t),u_t(t))\|_{X_1}^2.
\end{gather}

\noindent
{\bf Well-posedness and regularity results}
\\
Let $a\in(0,1).$ We recall that for any $(u^0,u^1)\in X_1,$ there exists a unique solution $(u,u_t)\in\vC([0,\infty);X_1)$ of (\ref{wave}). Moreover, $u(\: . \:, a)\in H^1_\loc([0,\infty)).$ Thus equation (\ref{wave}) makes sense in $L^2_\loc([0,\infty);H^{-1}(0,1)).$ In addition, $u$ satisfies the following energy estimate.
\begin{gather}
\label{deriveenergy}
\forall t\ge s\ge0, \; E_1(u(t))-E_1(u(s))=-\vint_s^t|u_t(\sigma,a)|^2d\sigma\le0.
\end{gather}
If furthermore $(u^0,u^1)\in\vD(\vA_1)$ then $(u,u_t)\in\vC([0,\infty);\vD(\vA_1)).$ Finally, $\vA_1$ is $m$--dissipative with domain dense in $X_1$ so that $\vA_1$ generates a semigroup of contractions $(\vS_1(t))_{t\ge0}$ on $X_1$ and on $\vD(\vA_1),$ which means that
\begin{gather}
\nonumber
\forall(u^0,u^1)\in X_1,\; \|(u(t),u_t(t))\|_{X_1}\le\|(u^0,u^1)\|_{X_1}, \\
\label{dec}
\forall(u^0,u^1)\in\vD(\vA_1),\; \|(u(t),u_t(t))\|_{\vD(\vA_1)}\le\|(u^0,u^1)\|_{\vD(\vA_1)},
\end{gather}
for any $t\ge0.$ For more details, see for example Theorem~1.1 and Lemma~2.1 of Tucsnak~\cite{MR1627724} and Proposition~2.1 of Ammari and Tucsnak~\cite{MR2002f:93104}. We also recall that $E_1(u(t))\xrightarrow{t\tends\infty}0,$ or equivalently
\begin{gather}
\nonumber
\vlim_{t\to\infty}\left(\|u(t)\|_{V_1} + \|u_t(t)\|_{L^2(0,1)}\right)=0, \\
\nonumber
\mbox{if and only if} \\
\label{tend0}
a\not\in\Q.
\end{gather}
And if furthermore $a$ satisfies (\ref{tend0}) and if $(u^0,u^1)\in\vD(\vA_1)$ then we have the estimate
$$
\forall t\ge0, \;
\|(u(t),u_t(t))\|_{X_1}\le\|\vS_1(t)\|_{\vL(\vD(\vA_1);X_1)}\|(u^0,u^1)\|_{\vD(\vA_1)},
$$
with $\vlim_{t\to\infty}\|\vS_1(t)\|_{\vL(\vD(\vA_1);X_1)}=0$ (Proposition~1.1 of Tucsnak~\cite{MR1627724}). Finally, it follows from (\ref{energy1})--(\ref{deriveenergy}) that
\begin{gather}
\label{decE}
\forall t\ge s\ge0, \;
\|(u(t),u_t(t))\|_{X_1}\le\|(u(s),u_t(s))\|_{X_1}.
\end{gather}
Our goal is to describe the decay rate of $E_1(u(t))$ as $t\tends\infty,$ for any $a\in(0,1)$ as soon as $E_1(u(t))\xrightarrow{t\tends\infty}0,$ when the lack of observability occurs. By (\ref{tend0}), this means that $a\not\in\Q.$
\medskip
\\
{\bf Known decay}

\noindent
Now, we show that our method allows us to recover the known results (Jaffard, Tucsnak and Zuazua~\cite{MR99g:93073}). We recall the definition of an irrational algebraic number.

\begin{vdefi}
Let $d\in\N,$ $d\ge2.$ An irrational number $a$ is said to be {\it algebraic of degree} $d$ if there exists a minimal polynomial function $P$ of degree $d$ with rational coefficients such that $P(a)=0.$ $P$ is {\it minimal} in the sense that if $Q$ is a polynomial function with rational coefficients such that $Q(a)=0$ then $\deg Q\ge\deg P.$
\end{vdefi}

\noindent
If $a$ is an irrational algebraic number of degree $d$ then it follows from Liouville's Theorem that there exists a positive constant $C=C(d)$ such that for any $(m,n)\in\Z\times\N,$ $\left|a-\frac{m}{n}\right|\ge\frac{C}{n^d}.$ This implies that there exists a positive constant $c_1=c_1(d)$ such that
\begin{gather}
\label{algebraic}
\forall n\in\N, \; |\sin(n\pi a)|\ge\frac{c_1}{n^{d-1}}
\quad \mbox{ and } \quad
\left|\sin\left(\left(n+\frac{1}{2}\right)\pi a\right)\right|\ge\frac{c_1}{(2n+1)^{d-1}}.
\end{gather}

\begin{vnota}
\label{S}
We denote by $\vS$ the set of all irrational numbers $a\in(0,1)$ such that if $[0,a_1,\ldots,a_n,\ldots]$ is the expansion of $a$ as a continued fraction, then $(a_n)_{n\in\N}$ is bounded.
\end{vnota}

\noindent
Let us notice that $\vS$ is obviously infinite and not countable and by classical results on Diophantine approximation (see Cassals~\cite{Cassals}, p.120), $\lambda(\vS)=0,$ where $\lambda$ is the Lebesgue's measure. Moreover, by Euler--Lagrange's Theorem (see Lang~\cite{MR35:129}, p.57), $\vS$ contains the set of algebraic irrational numbers $a\in(0,1)$ of degree 2. According to a classical result (see Tucsnak~\cite{MR1627724} and the references therein), if $a\in\vS$ then estimates (\ref{algebraic}) hold with $d=2.$ Finally, for any $\eps>0,$ there exist two $\lambda$--measurable sets $I_\eps\subset(0,1)$ and $J_\eps\subset(0,1)$ and a constant $c_2=c_2(\eps)>0$ such that $\lambda(I_\eps)=\lambda(J_\eps)=1$ and such that for any $a\in I_\eps$ and any $b\in J_\eps,$
\begin{gather}
\label{s}
\forall n\in\N, \; |\sin(n\pi a)|\ge\frac{c_2}{n^{1+\eps}}
\quad \mbox{ and } \quad
\left|\sin\left(\left(n+\frac{1}{2}\right)\pi b\right)\right|\ge\frac{c_2}{(2n+1)^{1+\eps}}.
\end{gather}
Let us notice that by Roth's Theorem (see Cassals~\cite{Cassals}, p.104), $I_\eps$ and $J_\eps$ contain all algebraic irrational numbers of $(0,1).$ The following result is due to Jaffard, Tucsnak and Zuazua~\cite{MR99g:93073} (Theorem~3.3).

\begin{vprop}[\cite{MR99g:93073}]
\label{propalgebraic}
Let $\vS$ be defined in Notation~$\ref{S}$ and let for any $t\ge0,$ $\omega_2(t)=t^2.$ We have the following result.
\begin{enumerate}
\item
\label{propalgbraic-1}
Let $a\in\vS$ and set for any $t>0,$ $\omega_1(t)=\frac{c_1}{t},$ where $c_1$ is given by $(\ref{algebraic})$ with $d=2.$ Then there exists a constant $C=C(a)>0$ such that for any initial data $(u^0,u^1)\in\vD(\vA_1),$ the corresponding solution $u$ of $(\ref{wave})$ verifies
\begin{gather}
\label{propalgebraic1}
E_1(u(t))\le\dfrac{C}{(t+1)}\|(u^0,u^1)\|_{\vD(\vA_1)}^2,
\end{gather}
for any $t\ge0.$ Furthermore, time decay in $(\ref{propalgebraic1})$ is optimal in the sense of Definition~$\ref{defioptim}.$
\item
\label{propalgbraic-2}
Let $\eps>0$ and set for any $t>0,$ $\omega_1(t)=\frac{c_2}{t^{1+\eps}},$ where $c_2$ is given by $(\ref{s}).$ For almost every $a\in(0,1)\cap\Q^{\mathrm c},$ there exists a constant $C=C(a,\eps)>0$ such that for any initial data $(u^0,u^1)\in\vD(\vA_1),$ the corresponding solution $u$ of $(\ref{wave})$ verifies
\begin{gather}
\label{propalgebraic2}
E_1(u(t))\le\dfrac{C}{(t+1)^\frac{1}{1+\eps}}\|(u^0,u^1)\|_{\vD(\vA_1)}^2,
\end{gather}
for any $t\ge0.$ Furthermore, time decay in $(\ref{propalgebraic2})$ is optimal in the sense of Definition~$\ref{defioptim}.$
\end{enumerate}
\end{vprop}

\subsection{New results}
\label{new}

Before stating the main results, let us make the following definition.

\begin{vdefi}
\label{adm}
We say that the functions $(\omega_1,\omega_2,\Phi,\Psi)$ are an {\it admissible quadruplet} if the following assertions hold.
\begin{enumerate}
\item
\label{adm1}
The quadruplet $(\omega_1,\omega_2,\Phi,\Psi)$ satisfies (\ref{hypo5})--(\ref{hypo8}) on $(0,\infty)$ and (\ref{hypo9}) holds on $(1,\infty).$
\item
\label{adm3}
One of the two following conditions is satisfied.
\begin{enumerate}
\item
\label{adm3-1}
The function $t\longmapsto\dfrac{1}{t}\vH_{\Phi,\Psi}^{-1}(t)$ is nondecreasing on $(0,1),$ where $\vH_{\Phi,\Psi}^{-1}$ defined by (\ref{prelem-2}) has to verify $\vH_{\Phi,\Psi}((0,\delta))\supset(0,1).$
\item
\label{adm3-2}
For any $t>0,$ $\Phi(t)=C_1t^\frac{1}{p}$ and $\Psi(t)=C_2t^\frac{1}{q}$ for some $p\ge1,$ $q\ge1$ and constants $C_1,C_2>0.$ In particular, we have for any $t>0,$ $\vH_{\Phi,\Psi}(t)=\left(C_1C_2^{-1}\right)^qt^\frac{q}{p}.$
\end{enumerate}
\end{enumerate}
\end{vdefi}

\noindent
In our applications, the weight $\omega_1$ comes from an oscillating function and it is not clear that it satisfies \eqref{hypo5}. So we precise how we obtain such a weight.

\begin{vlem}
\label{od}
Let $-\infty<a<b\le\infty$ and let $\eps:[a,b)\tends(0,\infty)$ be a continuous function such that
$\vliminf_{t\nearrow b}\eps(t)=0.$ Then there exists a convex function $\vphi\in\vC_\b^1([a,b);\R)$ such that $0<\vphi\le\eps$ and $\vphi^\prime<0$ on $[a,b).$
\end{vlem}

\begin{proof*}
Firstly, we note that we can find a positive function $\widetilde\eps\in\vC^1([a,b);\R)$ such that $0<\widetilde\eps\le\eps$ and $\widetilde\eps^{\:\prime}<0$ on $[a,b).$ So it is enough to consider $\eps$ to be such a function. Secondly, up to a bijective transformation conserving the convexity, we may assume that $[a,b)=[0,1).$ Set
\begin{gather*}
\forall t\in[0,1), \; f(t)=\max\{\eps^\prime(s);\; 0\le s\le t\}.
\end{gather*}
Define $\vphi$ by
\begin{gather*}
\forall t\in[0,1), \; \vphi(t)=-\vint_t^1f(s)ds \; \mbox{ and } \; \vphi(1)=0.
\end{gather*}
Since $f$ is monotone and $\eps^\prime$ is continuous, it follows that $f\in\vC_\b([0,1);\R).$ Then $\vphi$ is well-defined, $\vphi\in\vC_\b([0,1];\R)\cap\vC_\b^1([0,1);\R)$ and $\vphi^\prime=f$ on $[0,1).$ Clearly, $\vphi>0$ and $\vphi^\prime<0$ on $[0,1).$ In addition, $\vphi^\prime$ is nondecreasing so that $\vphi$ is convex. Finally, for any $\sigma\in[0,1),$ $\vphi^\prime(\sigma)\ge\eps^\prime(\sigma).$ Integrating this expression on $(t,1),$ for any $t\in[0,1),$ and using that $\vphi(1)=\eps(1)=0,$ we get $\vphi(t)\le\eps(t).$ This concludes the proof.
\medskip
\end{proof*}

\noindent
Let $(u_n)_{n\in\N}\subset(0,\infty)$ be such that $\vliminf_{n\to\infty}u_n=0.$ Let $\eps\in\vC([0,\infty);\R)$ be such that $0<\eps(n)\le u_n,$ for any $n\in\N.$ Let $\vphi\in\vC([0,\infty);\R)$ be a decreasing and convex function such that for any $t\ge0,$ $0<\vphi(t)\le\eps(t)$ (which exists by Lemma~\ref{od}) and consider $\vC\subset[1,\infty)\times[0,\infty)$ the closure of the convex envelope of the set $\{(n,u_n); \; n\in\N\}.$ Finally, fix arbitrarily $t\ge1.$ Then the set $\vC_t\stackrel{{\rm def}}{=}\vC\cap\left(\{t\}\times\R\right)$ is nonempty, closed and Lemma~\ref{od} ensures that for any $s_t\in\R$ such that $(t,s_t)\in\vC_t,$
\begin{gather*}
0<\vphi(t)\le s_t.
\end{gather*}
So by compactness, we may define the function $\omega_1$ as
\begin{gather}
\label{envc}
\forall t\ge1, \; \omega_1(t)=\min\{s_t;\; (t,s_t)\in\vC_t\}
\end{gather}
and extend $\omega_1$ as a decreasing, continuous and convex way on $[0,1].$ From the above discussion, Lemma~\ref{od} and Remark~\ref{rmkpropphi}, $\omega_1$ satisfies \eqref{hypo5} with $m=0.$ This justifies the following definition.

\begin{vdefi}
\label{o1}
Let $(u_n)_{n\in\N}\subset(0,\infty)$ such that $\vliminf_{n\to\infty}u_n=0.$ The function $\omega_1$ defined on $[0,\infty)$ by \eqref{envc} is called the {\it lower convex envelope of the sequence} $(u_n)_{n\in\N}.$
\end{vdefi}

\noindent
In some sense, $\omega_1$ is the ``nearest'' convex and decreasing function of $(u_n)_{n\in\N}$ satisfying $0<\omega_1(n)\le u_n,$ for any $n\in\N.$ It will be useful to consider the weights $\omega_1$ and $\omega_2$ defined as following. Let $a\in(0,1)\cap\Q^{\mathrm c}.$
\begin{gather}
\label{omega1}
\omega_1 \mbox{ is the lower convex envelope of the sequence } (\sin^2(n\pi a))_{n\in\N}, \medskip \\
\label{omega2}
\forall t\ge0,\; \omega_2(t)=t^2.
\end{gather}

\noindent
The following lemma shows that such definition for weights is consistent with the notion of admissible quadruplet.

\begin{vprop}
\label{propdefomega}
Let $(u_n)_{n\in\N}\subset(0,\infty)$ be such that $\vliminf_{n\to\infty}u_n=0,$ let $\omega_1$ be its lower convex envelope $($Definition~$\ref{od}),$ let $p\ge1,$ let $\alpha\in[0,1]$ and set for any $t\ge0,$ $\omega_2(t)=(t+\alpha)^p.$ Define for any $t\ge\alpha^p,$ $\Psi(t)=t^\frac{1}{p}-\alpha$ and for any $t>0,$
\begin{gather*}
\vphi(t)=\frac{\omega_1(t)}{t^p} \quad and \quad \Phi(t)=\frac{1}{\vphi^{-1}(t)}.
\end{gather*}
Then the quadruplet $(\omega_1,\omega_2,\Phi,\Psi)$ is admissible and for any $t>0,$
\begin{gather*}
\vH_{\Phi,\Psi}(t)=\frac{1}{\left(\vphi^{-1}(t)+\alpha\right)^p}.
\end{gather*}
\end{vprop}

\begin{proof*}
By definition of $\omega_1,$ $\omega_2$ and $\Psi,$ \eqref{hypo5}, \eqref{hypo6} and \eqref{hypo8} are satisfied. 
By Lemma~\ref{concave} applied to $f=\omega_1$ and with $m=0$ and $M=\omega(0),$ it follows that $\Phi$ satisfies \eqref{hypo7}. Moreover, we easily check that $\Phi\ge\frac{1}{\omega^{-1}}$ on $(0,\omega_1(1)].$ As a consequence, \eqref{hypo9} holds on $[1,\infty),$ so that condition~\ref{adm1} of Definition~\ref{adm} is fulfilled. Finally, by Lemma~\ref{prelem} we have
\begin{align*}
\forall t\in\left(0,\alpha^{-p}\right), & \;\; \widetilde\vH(t)\stackrel{{\rm def}}{=}
\frac{1}{t}\vH_{\Phi,\Psi}^{-1}(t)=\left(1-\alpha t^\frac{1}{p}\right)^{-p}\omega_1\left(t^{-\frac{1}{p}}-\alpha\right), \medskip \\
\forall t>0, & \;\;
\vH_{\Phi,\Psi}(t)=\frac{1}{\left(\vphi^{-1}(t)+\alpha\right)^p},
\end{align*}
where we used the notation $\alpha^{-p}=+\infty$ if $\alpha=0.$ It is clear that $\widetilde\vH$ is increasing on $\left(0,\alpha^{-p}\right)\supset(0,1),$ so that \eqref{adm3-1} of Definition~\ref{adm} holds and $(\omega_1,\omega_2,\Phi,\Psi)$ is an admissible quadruplet.
\medskip
\end{proof*}

\noindent
The main results are the following.

\begin{vthm}
\label{thmdecay1}
Let $a\in(0,1)\cap\Q^{\mathrm c}$ and let $\omega_1$ and $\omega_2$ be defined by
$(\ref{omega1})-(\ref{omega2}).$ Let $\Phi$ and $\Psi$ be two functions such that the quadruplet $(\omega_1,\omega_2,\Phi,\Psi)$ is admissible $($Definition~$\ref{adm}).$ Let $\vH_{\Phi,\Psi}$ be defined by $(\ref{optim1}).$ Then there exists a constant $C=C(a)>0$ such that for any initial data $(u^0,u^1)\in\vD(\vA_1),$ the corresponding solution $u$ of $(\ref{wave})$ verifies
\begin{gather}
\label{thmdecay0-1}
\forall t\ge0,\;
E_1(u(t))\le C\vH_{\Phi,\Psi}\left(\frac{1}{t+1}\right)\|(u^0,u^1)\|_{\vD(\vA_1)}^2,
\end{gather}
if $\Phi$ and $\Psi$ satisfy the hypothesis~$(\ref{adm3-1})$ of Definition~$\ref{adm}$ and
\begin{gather}
\label{thmdecay0-2}
\forall t\ge0,\;
E_1(u(t))\le\frac{C}{(t+1)^\frac{q}{p}}\|(u^0,u^1)\|_{\vD(\vA_1)}^2,
\end{gather}
if for any $t>0,$ $\Phi(t)=C_1t^\frac{1}{p}$ and $\Psi(t)=C_2t^\frac{1}{q}$ for some $p\in[1,\infty),$ $q\in[1,\infty)$ and constants $C_1,C_2>0$ $($case $(\ref{adm3-2})$ of Definition~$\ref{adm}).$
\end{vthm}

\begin{vrmk}
\label{rmkthmdecay1}
At the light of estimate (\ref{thmdecay0-1}), it is clear that we would like to find some functions $\Phi$ and $\Psi$ such that $\vH_{\Phi,\Psi}(t)$ goes to $0$ as $t\searrow0$ as rapidly as possible. This justifies Definition~\ref{defioptim}. Moreover, Proposition~\ref{propdefomega} ensures that there exists a quadruplet of functions $(\omega_1,\omega_2,\Phi,\Psi)$ which is admissible.
\end{vrmk}

\noindent
Concerning the explicit decay, the results are the following.

\begin{vthm}
\label{thmdecay2}
Let $a\in(0,1)\cap\Q^{\mathrm c}$ and let $\omega_1$ be defined by $(\ref{omega1}).$ We set
\begin{gather*}
\forall t>0, \; \vphi(t)=\dfrac{\omega_1(t)}{t^2}.
\end{gather*}
Then there exists a constant $C=C(a)>0$ such that for any initial data $(u^0,u^1)\in\vD(\vA_1),$ the corresponding solution $u$ of $(\ref{wave})$ satisfies
\begin{gather*}
\forall t\ge0, \;
\|(u(t),u_t(t))\|_{V_1\times L^2(0,1)}\le\dfrac{C}{\vphi^{-1}\left(\frac{1}{t+1}\right)}\|(u^0,u^1)\|_{\vD(\vA_1)}.
\end{gather*}
\end{vthm}

\begin{vrmk}
\label{rmkthmdecay2-1}
By Theorem~\ref{thmdecay2}, we are able to give the explicit decay of the energy for any
$a\in(0,1)\cap\Q^{\mathrm c}.$ This completes the lack, since the decay was known for almost every $a\in(0,1)$ (Jaffard, Tucsnak and Zuazua~\cite{MR99g:93073}, Theorem~3.3).
\end{vrmk}

\begin{vrmk}
\label{rmkthmdecay2-2}
It follows from Theorem~\ref{thmdecay2} and Proposition~\ref{propdefomega} that for any $(u^0,u^1)\in\vD(\vA_1),$ the corresponding solution $u$ of $(\ref{wave})$ satisfies
\begin{gather*}
\|(u(t),u_t(t))\|_{V_1\times L^2(0,1)}\le C\Phi\left(\frac{1}{t+1}\right)\|(u^0,u^1)\|_{\vD(\vA_1)}.
\end{gather*}
for any $t\ge0.$ In other words, decay of the energy directly depends on the behavior of the interpolation function $\Phi$ near $0.$
\end{vrmk}

\begin{vproof}{of Theorem~\ref{thmdecay2}.}
The result comes from Proposition~\ref{propdefomega} (applied with $(u_n)_{n\in\N}=(\sin^2(n\pi a))_{n\in\N},$ $p=2$ and $\alpha=0)$ and from (\ref{thmdecay0-1}) of Theorem~\ref{thmdecay1}.
\medskip
\end{vproof}

\begin{vproof}{of Proposition~\ref{propalgebraic}.}
Let $\vS$ be defined in Notation~\ref{S}. \\
{\bf Case of \ref{propalgbraic-1}.} Let $a\in\vS$ and let $c_1$ be the constant in (\ref{algebraic}) with $d=2.$ \\
{\bf Case of \ref{propalgbraic-2}.} Let $\eps>0,$ let $I_\eps\subset(0,1)$ be the set introduced after the Notation~\ref{S}, let $c_2$ be the constant in (\ref{s}) and let $a\in I_\eps.$ \\
{\bf Preliminary.}
Let $\nu\ge0$ and $\ell\in\{1,2\}.$ We define on $(0,\infty)$ the following functions.
\begin{gather*}
\omega_1(t)=\frac{c_\ell^2}{t^{2(1+\nu)}},
\quad
\Psi(t)=t^\frac{1}{2},
\quad
\Phi(t)=2\left(\frac{t}{c_\ell^2}\right)^\frac{1}{2(1+\nu)}.
\end{gather*}
Let $\omega_2$ be defined by (\ref{omega2}) and let $\vH_{\Phi,\Psi}$ be the corresponding functions given by (\ref{optim1}). Then
\begin{gather*}
\forall t>0,\; \vH_{\Phi,\Psi}(t)=4\left(\frac{t}{c_\ell^2}\right)^\frac{1}{1+\nu}.
\end{gather*}
Furthermore for any $t>0,$ $\Phi(\omega_1(t))\Psi(\omega_2(t))\ge1$ and $\vH_{\Phi,\Psi}(t)=\frac{C}{\omega_2\circ\omega_1^{-1}(t)}.$
\\
{\bf Proof of \ref{propalgbraic-1}.} 
Let $\nu=0$ and $\ell=1.$ The result follows by applying (\ref{thmdecay0-2}) of Theorem~\ref{thmdecay1}.
\\
{\bf Proof of \ref{propalgbraic-2}.} 
Let $\nu=\eps$ and $\ell=2.$ The result follows by applying (\ref{thmdecay0-2}) of Theorem~\ref{thmdecay1}. This concludes the proof.
\medskip
\end{vproof}

\noindent
Before proving Theorem~\ref{thmdecay1}, we need several results. Let us decompose the solution $u$ as following. For $u$ solution of (\ref{wave}) with initial data $(u^0,u^1)\in X_1,$ we write
\begin{gather}
\label{decompositionud}
u(t,x)=v(t,x)+w(t,x),
\end{gather}
for $(t,x)\in[0,\infty)\times(0,1),$ where $v$ is the unique solution of
\begin{gather}
 \left\{
  \begin{split}
   \label{wavev}
    v_{tt}-v_{xx}=0, & \quad (t,x)\in(0,\infty)\times(0,1), \\
      v(0,x)=u^0(x), & \quad  x\in(0,1), \\
    v_t(0,x)=u^1(x),& \quad  x\in(0,1), \\
     v(t,0)=v(t,1)=0, & \quad t\in[0,\infty).
  \end{split}
 \right.
\end{gather}
Then we have the well-known result (see for example Lemmas~4.1 and 5.3 of Ammari and Tucsnak~\cite{MR2002f:93104} for the proof).
\begin{vlem}
\label{equi}
Let $a\in(0,1)$ and let $T=10.$ Then there exists a constant $C_1=C_1(a)>0$ satisfying the following property. For any initial data $(u^0,u^1)\in X_1,$ the corresponding solutions $u$ and $v$ of $(\ref{wave})$ and $(\ref{wavev})$ satisfy
\begin{gather*}
C_1\vint_0^{\;T}v_t^2(t,a)dt\le\vint_0^{\;T}u_t^2(t,a)dt\le4\vint_0^{\;T}v_t^2(t,a)dt.
\end{gather*}
\end{vlem}

\noindent
Now, we decompose $u^0\in V_1$ and $u^1\in L^2(0,1)$ as
\begin{gather}
\label{fourieru0wd}
u^0(x)=\vsum_{n=0}^\infty a_n\sin(n\pi x), \quad
u^1(x)=\pi\vsum_{n=0}^\infty nb_n\sin(n\pi x).
\end{gather}
We then have
\begin{gather}
\label{fouriernormwd}
\|u^0\|_{L^2(0,1)}^2=\frac{1}{2}\vsum_{n=0}^\infty a_n^2, \quad
\|u_x^0\|_{L^2(0,1)}^2=\frac{\pi^2}{2}\vsum_{n=0}^\infty n^2a_n^2, \quad
\|u^1\|_{L^2(0,1)}^2=\frac{\pi^2}{2}\vsum_{n=0}^\infty n^2b_n^2.
\end{gather}
It follows that the solution $v$ of (\ref{wavev}) is defined by
\begin{gather}
\label{fouriervwd}
\forall(t,x)\in\R\times(0,1), \;
v(t,x)=\vsum_{n=0}^\infty\left\{(a_n\cos(n\pi t) + b_n \sin(n\pi t))\sin(n\pi x)\right\}.
\end{gather}
If furthermore $(u^0,u^1)\in D(A_1)\times V_1$ then
\begin{gather}
\label{fouriernormdregularwd}
\|u^0_{xx}\|_{L^2(0,a)}^2+\|u^0_{xx}\|_{L^2(a,1)}^2=\frac{\pi^4}{2}\vsum_{n=0}^\infty n^4a_n^2,
\quad
\|u^1_x\|_{L^2(0,1)}^2=\frac{\pi^4}{2}\vsum_{n=0}^\infty n^4b_n^2.
\end{gather}

\noindent
We have the following simple result.

\begin{vlem}
\label{borneinfv}
Let $a\in(0,1),$ let $T=10,$ let $(u^0,u^1)\in X_1$ and let $(a_n)_{n\in\N}\in\ell^2(\N)$ and $(b_n)_{n\in\N}\in\ell^2(\N)$ be given by $(\ref{fourieru0wd}).$ Then
\begin{gather}
\label{borninfv1}
\vint_0^{\;T}v_t^2(t,a)dt\ge\pi^2\vsum_{n=0}^\infty n^2(a_n^2+b_n^2)\sin^2(n\pi a),
\end{gather}
where $v$ is the solution of $(\ref{wavev})$ given by $(\ref{fouriervwd}).$
\end{vlem}

\begin{proof*}
Using (\ref{fouriervwd}), we have
\begin{align*}
          \vint_0^{\;\;T}v_t^2(t,a)dt
\ge \; & \pi^2\vint_0^2\left(\vsum_{n=0}^\infty\sin(n\pi a)(-na_n\sin(n\pi t)+nb_n\cos(n\pi t))\right)^2dt \medskip \\
  =  \; & \pi^2\vsum_{n=0}^\infty\sin^2(n\pi a)(n^2a_n^2+n^2b_n^2),
\end{align*}
where the last line comes from Parseval's identity. Hence (\ref{borninfv1}).
\medskip
\end{proof*}

\begin{vlem}
\label{estiminf}
Let $a\in(0,1)\cap\Q^{\mathrm c},$ let $T=10,$ let $\omega_1$ be given by $(\ref{omega1})$ and let $\omega_2$ be given by $(\ref{omega2}).$ Let $\Phi$ and $\Psi$ be two functions such that the quadruplet $(\omega_1,\omega_2,\Phi,\Psi)$ satisfies hypothesis $\ref{adm1}$ of Definition~$\ref{adm}$ and such that $\vH_{\Phi,\Psi}((0,\delta))\supset(0,1).$ Then there exists a constant $C_2=C_2(a)>0$ such that for any initial data $(u^0,u^1)\in\vD(\vA_1),$
\begin{gather}
\label{estiminf1}
\|(u^0,u^1)\|_{X_1}^2-\|(u(T),u_t(T))\|_{X_1}^2\ge C_2\|(u^0,u^1)\|_{X_1}^2\vH_{\Phi,\Psi}^{-1}
\left(\frac{\|(u^0,u^1)\|_{X_1}^2}{\|(u^0,u^1)\|_{\vD(\vA_1)}^2}\right),
\end{gather}
where $u$ is the solution of $(\ref{wave}),$ and where $\vH_{\Phi,\Psi}^{-1}$ is defined by $(\ref{prelem-2}).$
\end{vlem}

\begin{proof*}
By Proposition~\ref{propdefomega}, $\Phi$ and $\Psi$ exist. We decompose $u^0$ and $u^1$ as in (\ref{fourieru0wd}). We write
\begin{gather}
\label{energy-}
E_-(0)=\dfrac{\pi^2}{2}\dsp\vsum_{n=0}^\infty n^2(a_n^2+b_n^2)\omega_1(n),
\end{gather}
where $\omega_1$ verifies \ref{adm1} of Definition~\ref{adm}. By (\ref{deriveenergy}) and Lemmas~\ref{equi} and \ref{borneinfv}, there exists a constant $C_2=C_2(a)>0$ such that
\begin{gather}
\label{estimateenergy-v}
\|(u^0,u^1)\|_{X_1}^2-\|(u(T),u_t(T))\|_{X_1}^2\ge C_2E_-(0).
\end{gather}
Assume further that $(u^0,u^1)\in D(A_1)\times V_1.$ We define
\begin{gather}
\label{energytemp+}
E_+(0)=\dfrac{\pi^4}{4}\dsp\vsum_{n=0}^\infty n^4(a_n^2+b_n^2).
\end{gather}
Putting together (\ref{energytemp+}) and (\ref{fouriernormdregularwd}), we have that for any initial data $(u^0,u^1)\in D(A_1)\times V_1,$
\begin{gather*}
E_+(0)=
\dfrac{1}{2}\left(\|u^0_{xx}\|_{L^2(0,a)}^2+\|u^0_{xx}\|_{L^2(a,1)}^2+\|u^1_x\|_{L^2(0,1)}^2\right).
\end{gather*}
These estimates imply that
\begin{gather}
\label{energy+}
E_+(0)\le\|(u^0,u^1)\|_{\vD(\vA_1)}^2.
\end{gather}
Recall that by (\ref{energy1}) and (\ref{fouriernormwd}),
\begin{gather}
\label{energy}
E(0)=\dfrac{\pi^2}{4}\dsp\vsum_{n=0}^\infty n^2(a_n^2+b_n^2)
\stackrel{{\rm def}}{=}\frac{1}{2}\|(u^0,u^1)\|_{X_1}^2,
\end{gather}
where we have set $E(0)=E_1(u(0)).$ Let $u=(u_n)_{n\in\N}\in\ell^1(\N;\R)$ be defined by
\begin{gather*}
\forall n\in\N, \; u_n=n^2(a_n^2+b_n^2).
\end{gather*}
Then it follows from Theorem~\ref{ineg1} (applied to the function $f=u,$ with $p=1,$ the discrete measure on $\vP(\N)$ and the weights $\omega_1$ and $\omega_2),$ (\ref{energy-}) and (\ref{energytemp+})--(\ref{energy}) that
$$
1\le\Phi\left(\frac{E_-(0)}{\|(u^0,u^1)\|_{X_1}^2}\right)
\Psi\left(\frac{\|(u^0,u^1)\|_{\vD(\vA_1)}^2}{\|(u^0,u^1)\|_{X_1}^2}\right),
$$
which yields
$$
E_-(0)\ge\|(u^0,u^1)\|_{X_1}^2\Phi^{-1}
\left(\dfrac{1}{\Psi\left(\dfrac{\|(u^0,u^1)\|_{\vD(\vA_1)}^2}{\|(u^0,u^1)\|_{X_1}^2}\right)}\right).
$$
Then for any $(u^0,u^1)\in D(A_1)\times V_1,$ 
\begin{gather}
\label{proofestiminf}
E_-(0)\ge\|(u^0,u^1)\|_{X_1}^2\vH_{\Phi,\Psi}^{-1}
\left(\frac{\|(u^0,u^1)\|_{X_1}^2}{\|(u^0,u^1)\|_{\vD(\vA_1)}^2}\right).
\end{gather}
From (\ref{estimateenergy-v}) and (\ref{proofestiminf}), it follows that (\ref{estiminf1}) holds for any $(u^0,u^1)\in D(A_1)\times V_1.$ By continuity of $\vH_{\Phi,\Psi}^{-1}$ and by density of $D(A_1)\times V_1$ in $Y_1$ (which contains $\vD(\vA_1)$ and has the same norm of $\vD(\vA_1)),$ it follows that (\ref{estiminf1}) holds for any $(u^0,u^1)\in\vD(\vA_1).$ Hence the result.
\medskip
\end{proof*}

\begin{vproof}{of Theorem~\ref{thmdecay1}.}
We follow the proof of Theorem~2.4 of Ammari and Tucsnak~\cite{MR2002f:93104}. Let $T=10.$ By Lemma~\ref{estiminf}, we have that
$$
\|(u(T),u_t(T))\|_{X_1}^2\le\|(u^0,u^1)\|_{X_1}^2
-C_2\|(u^0,u^1)\|_{X_1}^2\vH_{\Phi,\Psi}^{-1}
\left(\frac{\|(u^0,u^1)\|_{X_1}^2}{\|(u^0,u^1)\|_{\vD(\vA_1)}^2}\right).
$$
This estimate remains valid in successive intervals $[\ell T,(\ell+1)T].$ So with (\ref{dec}), (\ref{decE}) and the fact that $\vH_{\Phi,\Psi}^{-1}$ is increasing (Lemma~\ref{prelem}), we obtain that
\begin{gather}
\begin{split}
\label{proofthmdecay1-1}
\|(u((\ell+1)T),u_t((\ell+1)T))\|_{X_1}^2\le & \|(u(\ell T),u_t(\ell T))\|_{X_1}^2 \medskip \\
-C_2\|(u(\ell & T),u_t(\ell T))\|_{X_1}^2\vH_{\Phi,\Psi}^{-1}
\left(\frac{\|(u((\ell+1)T),u_t((\ell+1)T))\|_{X_1}^2}{\|(u^0,u^1)\|_{\vD(\vA_1)}^2}\right),
\end{split}
\end{gather}
for every $\ell\in\N\cup\{0\}.$ \medskip \\
{\bf Case 1.} The functions $\Phi$ and $\Psi$ satisfy hypothesis (\ref{adm3-1}) of Definition~\ref{adm}. \\
Our expression (\ref{proofthmdecay1-1}) is the same that (4.16) in Ammari and Tucsnak~\cite{MR2002f:93104} (with $X\times V=X_1,$ $\|\: . \:\|_{Y_1\times Y_2}=\|\: . \:\|_{\vD(\vA_1)},$ $\vG=\vH_{\Phi,\Psi}^{-1}$ and $\theta=\frac{1}{2}).$ The rest of the proof follows as in~\cite{MR2002f:93104} (where (\ref{adm3-1}) of Definition~\ref{adm} is used). Then (\ref{thmdecay0-1}) follows. \medskip \\
{\bf Case 2.} The functions $\Phi$ and $\Psi$ satisfy hypothesis (\ref{adm3-2}) of Definition~\ref{adm}. \\
It follows that for any $t>0,$ $\vH_{\Phi,\Psi}^{-1}(t)=C_3t^\frac{p}{q}.$ Using again (\ref{decE}) and the definition of $\vH_{\Phi,\Psi}^{-1},$ (\ref{proofthmdecay1-1}) becomes
\begin{gather}
\begin{split}
\label{proofthmdecay1-2}
\|(u((\ell+1)T),u_t((\ell+1)T))\|_{X_1}^2\le\|(u(\ell T),u_t(\ell T)) & \|_{X_1}^2 \medskip \\
-C_4 & \frac{\|(u((\ell+1)T),u_t((\ell+1)T))\|_{X_1}^{2\frac{p+q}{q}}}
{\|(u^0,u^1\|_{\vD(\vA_1)}^{2\frac{p}{q}}},
\end{split}
\end{gather}
for every $\ell\in\N\cup\{0\}.$ Our expression (\ref{proofthmdecay1-2}) is the same that (4.23) in Ammari and Tucsnak~\cite{MR2002f:93104} (with $X\times V=X_1,$ $\|\: . \:\|_{Y_1\times Y_2}=\|\: . \:\|_{\vD(\vA_1)}$ and $\theta=\frac{q}{p+q}).$ The rest of the proof follows as in~\cite{MR2002f:93104}.
\medskip
\end{vproof}

\begin{vrmk}
We are not able to apply directly Theorem~2.4 of Ammari and Tucsnak~\cite{MR2002f:93104}. Indeed, in their theorem, the assumption (2.8) is
$$
\vint_0^2v_t^2(t,a)dt\ge C\|(u^0,u^1)\|_{V_1\times L^2(0,1)}^2\vG
\left(\frac{\|(u^0,u^1)\|_{L^2(0,1)\times H^{-1}(0,1)}^2}{\|(u^0,u^1)\|_{V_1\times L^2(0,1)}^2}\right),
$$
(where $\vG=\vH_{\Phi,\Psi}^{-1})$ and we can only show the weaker estimate (by the inequalities of interpolation)
$$
\vint_0^2v_t^2(t,a)dt\ge C\|(u^0,u^1)\|_{V_1\times L^2(0,1)}^2\vG
\left(\frac{\|(u^0,u^1)\|_{V_1\times L^2(0,1)}^2}{\|(u^0,u^1)\|_{\vD(\vA_1)}^2}\right).
$$
\end{vrmk}

\section{Others applications}
\label{damping2}

\subsection{Wave equation with mixed boundary condition}
\label{sswaves}

We consider a wave equation with a damping control concentrated on an interior point $a\in(0,1)$ with a homogenous Dirichlet boundary condition at the left end and a homogenous Neumann boundary condition at the right end,
\begin{gather}
 \left\{
  \begin{split}
   \label{waves}
    u_{tt}-u_{xx}+\delta_{a}u_t(t,a)=0, & \quad (t,x)\in(0,\infty)\times(0,1), \\
      u(0,x)=u^0(x), \; u_t(0,x)=u^1(x), & \quad x\in(0,1), \\
                                u(t,0)=u_x(t,1)=0, & \quad t\in[0,\infty).
   \end{split}
 \right.
\end{gather}

\noindent
{\bf Notations for the wave equation (\ref{waves}) with homogenous mixed Dirichlet and Neumann boundary condition} \\
Let $V_2=\left\{u\in H^1(0,1);\; u(0)=0\right\}.$ A direct calculation gives that for any $u\in V_2,$ $\|u\|_{L^2(0,1)}\le\|u_x\|_{L^2(0,1)},$ so we may endow $V_2$ of the norm $\|u\|_{V_2}=\|u_x\|_{L^2(0,1)},$ for any $u\in V_2.$ Let $X_2=V_2\times L^2(0,1),$
\begin{gather*}
Y_2=\left\{u\in V_2\cap H^2(0,a)\cap H^2(a,1);\; \frac{{\mathrm d}u}{{\mathrm d}x}(1)=0\right\}\times V_2,
\medskip \\
D(A_2)=\left\{u\in V_2\cap H^2(0,1);\; \frac{{\mathrm d}u}{{\mathrm d}x}(1)=0\right\},
\quad
A_2=-\frac{{\mathrm d}^2}{{\mathrm d}x^2}, \medskip \\
\vD(\vA_2)=\left\{(u,v)\in Y_2;\; \frac{{\mathrm d}u}{{\mathrm d}x}(a_+)
-\frac{{\mathrm d}u}{{\mathrm d}x}(a_-)=v(a)\right\},
\end{gather*}
with
$$
\|(u,v)\|_{\vD(\vA_2)}^2=\|(u,v)\|_{Y_2}^2=\|u\|_{H^2(0,a)}^2+\|u\|_{H^2(a,1)}^2+\|v\|_{H^1(0,1)}^2,
$$
and let
$
\begin{array}{rl}
\vA_2= & \!\!\!\! \left(
 \begin{array}{cc}
    0    &  \rm{Id}    \\
-A_2  &  -\delta_{a}
 \end{array}
\right).
\end{array}
$
We define the energy $E_2$ for $u$ solution of equation (\ref{waves}) by (\ref{energy1}).
\medskip
\\
{\bf Well-posedness and regularity results}
\\
Let $a\in(0,1).$ We recall that for any $(u^0,u^1)\in X_2,$ there exists a unique solution $(u,u_t)\in\vC([0,\infty);X_2)$ of (\ref{waves}). Moreover, $u(\: . \:, a)\in H^1_\loc([0,\infty)).$ Thus equation (\ref{waves}) makes sense in $L^2_\loc([0,\infty);H^{-1}(0,1)).$ In addition, $u$ satisfies the following energy estimate.
\begin{gather}
\label{deriveenergys}
\forall t\ge s\ge0, \; E_2(u(t))-E_2(u(s))=-\vint_s^t|u_t(\sigma,a)|^2d\sigma\le0.
\end{gather}
If furthermore $(u^0,u^1)\in\vD(\vA_2)$ then $(u,u_t)\in\vC([0,\infty);\vD(\vA_2)).$ Finally, $\vA_2$ is $m$--dissipative with domain dense in $X_2$ so that $\vA_2$ generates a semigroup of contractions $(\vS_2(t))_{t\ge0}$ on $X_2$ and on $\vD(\vA_2),$ which means that
\begin{gather*}
\forall(u^0,u^1)\in X_2,\; \|(u(t),u_t(t))\|_{X_2}\le\|(u^0,u^1)\|_{X_2}, \\
\forall(u^0,u^1)\in\vD(\vA_2),\; \|(u(t),u_t(t))\|_{\vD(\vA_2)}\le\|(u^0,u^1)\|_{\vD(\vA_2)},
\end{gather*}
for any $t\ge0.$ For more details, see Proposition~1.1 and Section~3 p.223 of Ammari, Henrot and Tucsnak~\cite{MR2002j:93073}. We also recall that $E_2(u(t))\xrightarrow{t\tends\infty}0,$ or equivalently
\begin{gather}
\nonumber
\vlim_{t\to\infty}\left(\|u(t)\|_{V_2} + \|u_t(t)\|_{L^2(0,1)}\right)=0 \\
\nonumber
\mbox{if and only if} \\
\label{tend0s}
\forall (p,q)\in\N\times\N, \; a\not=\dfrac{2p}{2q-1},
\end{gather}
And if furthermore $a$ satisfies (\ref{tend0s}) and if $(u^0,u^1)\in\vD(\vA_2)$ then we have the estimate
$$
\forall t\ge0, \;
\|(u(t),u_t(t))\|_{X_2}\le\|\vS_2(t)\|_{\vL(\vD(\vA_2);X_2)}\|(u^0,u^1)\|_{\vD(\vA_2)},
$$
with $\vlim_{t\to\infty}\|\vS_2(t)\|_{\vL(\vD(\vA_2);X_2)}=0$ (Proposition~3.1 of Ammari, Henrot and Tucsnak~\cite{MR2002j:93073}). Finally,
\begin{gather}
\begin{cases}
\nonumber
\exists\omega>0,\; \exists C=C(\omega)>0\; \mbox{ such that }\; \forall(u^0,u^1)\in X_2, \medskip \\
\forall t\ge0,\; E_2(u(t))\le Ce^{-\omega t}E_2(u(0))
\medskip 
\end{cases}
\\
\nonumber
\mbox{if and only if}
\medskip \\
\label{exps}
a=\dfrac{2p-1}{q},\; \mbox{ for some }\; (p,q)\in\N\times\N.
\end{gather}
See Theorem~1.2 of Ammari, Henrot and Tucsnak~\cite{MR2002j:93073}. It follows from (\ref{energy1}) and (\ref{deriveenergys}) that
\begin{gather*}
\forall t\ge s\ge0, \;
\|(u(t),u_t(t))\|_{X_2}\le\|(u(s),u_t(s))\|_{X_2}.
\end{gather*}
We are concerned by the decay rate of the energy $E_2(u(t))$ when it is not exponentially stable. In particular, by \eqref{tend0s} and \eqref{exps} this implies that $a\not\in\Q.$
\medskip
\\
The main results are the following.

\begin{vthm}
\label{thmdecay1s}
Let $a\in(0,1)\cap\Q^{\mathrm c}$ and let $\omega_1$ be the lower convex envelope of the sequence
$$
\left(\sin^2\left(\left(n+\frac{1}{2}\right)\pi a\right)\right)_{n\in\N}
$$
$($Definition~$\ref{o1}).$ Let $\omega_2$ be defined on $[0,\infty)$ by $\omega_2(t)=\left(t+\frac{1}{2}\right)^2.$ Let $\Phi$ and $\Psi$ be two functions such that the quadruplet $(\omega_1,\omega_2,\Phi,\Psi)$ is admissible $($Definition~$\ref{adm}).$ Let $\vH_{\Phi,\Psi}$ be defined by $(\ref{optim1}).$ Then there exists a constant $C=C(a)>0$ such that for any initial data $(u^0,u^1)\in\vD(\vA_2),$ the corresponding solution $u$ of $(\ref{waves})$ verifies
\begin{gather*}
\forall t\ge0,\;
E_2(u(t))\le C\vH_{\Phi,\Psi}\left(\frac{1}{t+1}\right)\|(u^0,u^1)\|_{\vD(\vA_2)}^2,
\end{gather*}
if $\Phi$ and $\Psi$ satisfy the hypothesis~$(\ref{adm3-1})$ of Definition~$\ref{adm}$ and
\begin{gather}
\label{thmdecay0-2s}
\forall t\ge0,\;
E_2(u(t))\le\frac{C}{(t+1)^\frac{q}{p}}\|(u^0,u^1)\|_{\vD(\vA_2)}^2,
\end{gather}
if for any $t>0,$ $\Phi(t)=C_1t^\frac{1}{p}$ and $\Psi(t)=C_2t^\frac{1}{q}$ for some $p\in[1,\infty),$ $q\in[1,\infty)$ and constants $C_1,C_2>0$ $($case $(\ref{adm3-2})$ of Definition~$\ref{adm}).$
\end{vthm}

\begin{proof*}
We write $u^0(x)=\vsum_{n=0}^\infty a_n\sin\left(\left(n+\frac{1}{2}\right)\pi x\right)$ and
$u^1(x)=\pi\vsum_{n=0}^\infty\left(n+\frac{1}{2}\right)b_n\sin\left(\left(n+\frac{1}{2}\right)\pi x\right)$ and we consider the solution $v$ of (\ref{wavev}) satisfying the same boundary condition as $u.$ We follow the method as for (\ref{wave}). Then from Ingham's inequality (Ingham~\cite{MR1545625}) and the results of Ammari, Henrot and Tucsnak~\cite{MR2002j:93073} (Lemma~4.2 of~\cite{MR2002j:93073}; see also Lemma~2.5 of~\cite{MR2002j:93073} and Lemma~4.1 of~\cite{MR2002f:93104}), we obtain for $T=10,$ 
\begin{gather*}
\vint_0^{\;T}u_t^2(t,a)dt\ge C(a)\vint_0^{\;T}v_t^2(t,a)dt
\ge C(a)\pi^2\vsum_{n=0}^\infty\left(n+\frac{1}{2}\right)^2(a_n^2+b_n^2)
\sin^2\left(\left((n+\frac{1}{2}\right)\pi a\right).
\end{gather*}
Then we define
\begin{gather*}
E_+(0)=\dfrac{\pi^4}{4}\dsp\vsum_{n=0}^\infty\left(n+\frac{1}{2}\right)^4(a_n^2+b_n^2)
=\dfrac{\pi^4}{4}\dsp\vsum_{n=0}^\infty\left(n+\frac{1}{2}\right)^2(a_n^2+b_n^2)\omega_2(n), \medskip \\
E(0)=\dfrac{\pi^2}{4}\dsp\vsum_{n=0}^\infty\left(n+\frac{1}{2}\right)^2(a_n^2+b_n^2)
\stackrel{{\rm def}}{=}\frac{1}{2}\|(u^0,u^1)\|_{X_2}^2, \medskip \\
E_-(0)=\dfrac{\pi^2}{2}\dsp\vsum_{n=0}^\infty\left(n+\frac{1}{2}\right)^2(a_n^2+b_n^2)\omega_1(n).
\end{gather*}
The result follows from the discussion at the beginning of Section~\ref{damping}.
\medskip
\end{proof*}

\noindent
Using Theorem~\ref{thmdecay1s} and Proposition~\ref{propdefomega} (applied with $(u_n)_{n\in\N}=\left(\sin^2\left(\left(n+\frac{1}{2}\right)\pi a\right)\right)_{n\in\N},$ $p=2$ and $\left.\alpha=\frac{1}{2}\right),$ we obtain the following result.

\begin{vthm}
\label{thmdecay2s}
Let $a\in(0,1)\cap\Q^{\mathrm c}$ and let $\omega_1$ and $\omega_2$ be defined as in Theorem~$\ref{thmdecay1s}.$ We set
\begin{gather*}
\forall t>0, \; \vphi(t)=\dfrac{\omega_1(t)}{t^2}.
\end{gather*}
Then there exists a constant $C=C(a)>0$ such that for any initial data $(u^0,u^1)\in\vD(\vA_2),$ the corresponding solution $u$ of $(\ref{waves})$ satisfies
\begin{gather*}
\forall t\ge0, \;
\|(u(t),u_t(t))\|_{V_2\times L^2(0,1)}\le\dfrac{C}{\vphi^{-1}\left(\frac{1}{t+1}\right)} \|(u^0,u^1)\|_{\vD(\vA_2)}.
\end{gather*}
\end{vthm}

\begin{vrmk}
By Theorem~\ref{thmdecay2s}, we are able to give the explicit decay of the energy for any $a\in(0,1)\cap\Q^{\mathrm c}.$ This completes the lack, since the decay was known for almost every $a\in(0,1),$ as stated in Theorem~1.4 of ~Ammari, Henrot and Tucsnak~\cite{MR2002j:93073}. In addition, with help of  \eqref{thmdecay0-2s} of Theorem~\ref{thmdecay1s}, our method allows us to recover the results of that Theorem~1.4.
\end{vrmk}

\subsection{Bernoulli--Euler beam with a pointwise interior damping control}
\label{ssbeam}

We consider a Bernoulli--Euler beam with a damping control concentrated in an interior point $a\in(0,1),$
\begin{gather}
 \left\{
  \begin{split}
   \label{beam}
      u_{tt}+u_{xxxx}+\delta_{a}u_t(t,a)=0, & \quad (t,x)\in(0,\infty)\times(0,1), \\
      u(0,x)=u^0(x), \quad u_t(0,x)=u^1(x), & \quad x\in(0,1), \\
    u(t,0)=u(t,1)=u_{xx}(t,0)=u_{xx}(t,1)=0, & \quad t\in[0,\infty).
  \end{split}
 \right.
\end{gather}
We also could have chosen the boundary condition
$$
\forall t\ge0, \;
u(t,0)=u_x(t,1)=u_{xx}(t,0)=u_{xxx}(t,1)=0,
$$
as in Ammari and Tucsnak~\cite{MR1814271}. But for conciseness of the paper, we do not consider this case.
\medskip
\\
{\bf Notations for the Bernoulli--Euler beam equation (\ref{beam})} \\
Let $V_3=H^1_0(0,1)\cap H^2(0,1).$ By Cauchy--Schwarz's inequality, we have $\|u\|_{L^2(0,1)}\le\|u_x\|_{L^2(0,1)}\le\|u_{xx}\|_{L^2(0,1)},$ for any $u\in V_3.$ So we may endow $V_3$ of the norm $\|u\|_{V_3}=\|u_{xx}\|_{L^2(0,1)},$ for any $u\in V_3.$ Let $X_3=V_3\times L^2(0,1),$
\begin{gather*}
Y_3=\left\{u\in H^1_0(0,1)\cap H^2(0,1)\cap H^4(0,a)\cap H^4(a,1);\; \frac{{\mathrm d}^2u}{{\mathrm d}x^2}(0)=\frac{{\mathrm d}^2u}{{\mathrm d}x^2}(1)=0\right\}\times V_3,
\medskip \\
D(A_3)=\left\{u\in H^1_0(0,1)\cap H^4(0,1);\; \frac{{\mathrm d}^2u}{{\mathrm d}x^2}(0)
=\frac{{\mathrm d}^2u}{{\mathrm d}x^2}(1)=0\right\}, \quad A_3=\frac{{\mathrm d}^4}{{\mathrm d}x^4},
\medskip \\
\vD(\vA_3)=\left\{(u,v)\in Y_3;\; \frac{{\mathrm d}^2u}{{\mathrm d}x^2}(a_+)=
\frac{{\mathrm d}^2u}{{\mathrm d}x^2}(a_-)\; \mbox{ and }\; \frac{{\mathrm d}^3u}{{\mathrm d}x^3}(a_+)
-\frac{{\mathrm d}^3u}{{\mathrm d}x^3}(a_-)=-v(a)\right\},
\end{gather*}
with
$$
\|(u,v)\|_{\vD(\vA_3)}^2=\|(u,v)\|_{Y_3}^2=\|u\|_{H^4(0,a)}^2+\|u\|_{H^4(a,1)}^2+\|v\|_{H^2(0,1)}^2,
$$
and let
$
\begin{array}{rl}
\vA_3= & \!\!\!\! \left(
 \begin{array}{cc}
     0    &  \rm{Id}    \\
 -A_3  &  -\delta_{a}
 \end{array}
\right).
\end{array}
$
We define the energy $E_3$ for $u$ solution of equation (\ref{beam}) by
\begin{gather}
\label{energy3}
\forall t\ge0, \; E_3(u(t))=\frac{1}{2}\left(\|u_t(t)\|_{L^2(0,1)}^2 + \|u_{xx}(t)\|_{L^2(0,1)}^2\right)
=\frac{1}{2}\|(u(t),u_t(t))\|_{X_3}^2.
\end{gather}

\noindent
{\bf Well-posedness and regularity results}
\\
We recall that for any $(u^0,u^1)\in X_3,$ there exists a unique solution $(u,u_t)\in\vC([0,\infty);X_3)$ of (\ref{beam}). Moreover, $u(\: . \:, a)\in H^1_\loc([0,\infty))$ and thus equation (\ref{beam}) makes sense in $L^2_\loc([0,\infty);H^{-2}).$ In addition, $u$ satisfies the following energy estimate.
\begin{gather}
\label{deriveenergy3}
\forall t\ge s\ge0, \; E_3(u(t))-E_3(u(s))=-\vint_s^t|u_t(\sigma,a)|^2d\sigma\le0.
\end{gather}
If furthermore $(u^0,u^1)\in\vD(\vA_3)$ then $(u,u_t)\in\vC([0,\infty);\vD(\vA_3)).$ Finally, $\vA_3$ is $m$--dissipative with domain dense in $X_3$ so that $\vA_3$ generates a semigroup of contractions $(\vS_3(t))_{t\ge0}$ on $X_3$ and on $\vD(\vA_3),$ which means that
\begin{gather*}
\forall(u^0,u^1)\in X_3,\; \|(u(t),u_t(t))\|_{X_3}\le\|(u^0,u^1)\|_{X_3}, \\
\forall(u^0,u^1)\in\vD(\vA_3),\; \|(u(t),u_t(t))\|_{\vD(\vA_3)}\le\|(u^0,u^1)\|_{\vD(\vA_3)},
\end{gather*}
for any $t\ge0.$ For more details, see for example Proposition~2.1 of Ammari and Tucsnak~\cite{MR2002f:93104}; Section~2 p.1161, Proposition~2.1 and Section~5 p.1173--1174 of Ammari and Tucsnak~\cite{MR1814271}. We also recall that $E_3(u(t))\xrightarrow{t\tends\infty}0,$ or equivalently
\begin{gather}
\nonumber
\vlim_{t\to\infty}\left(\|u(t)\|_{V_3} + \|u_t(t)\|_{L^2(0,1)}\right)=0 \\
\nonumber
\mbox{if and only if} \\
\label{tend03}
a\not\in\Q.
\end{gather}
And if furthermore $a$ satisfies (\ref{tend03}) and if $(u^0,u^1)\in\vD(\vA_3)$ then we have the estimate
\begin{gather*}
\forall t\ge0, \;
\|(u(t),u_t(t))\|_{X_3}\le\|\vS_3(t)\|_{\vL(\vD(\vA_3);X_3)}\|(u^0,u^1)\|_{\vD(\vA_3)},
\end{gather*}
with $\vlim_{t\to\infty}\|\vS_3(t)\|_{\vL(\vD(\vA_3);X_3)}=0$ (Proposition~2.1 and Section~5 p.1174 of Ammari and Tucsnak~\cite{MR1814271}). Finally, it follows from (\ref{energy3})--(\ref{deriveenergy3}) that the following holds.
\begin{gather*}
\forall t\ge s\ge0, \;
\|(u(t),u_t(t))\|_{X_3}\le\|(u(s),u_t(s))\|_{X_3}.
\end{gather*}
The goal is to establish the decay rate of $E_3(u(t))$ as $t\tends\infty,$ for any $a\in(0,1)$ as soon as $E_3(u(t))\xrightarrow{t\tends\infty}0,$ when the lack of observability occurs. In particular, by (\ref{tend03}), this implies that $a\not\in\Q.$

\begin{vthm}
\label{thmdecay1-2}
Let $a\in(0,1)\cap\Q^{\mathrm c},$ let $\omega_1$ be the lower convex envelope of the sequence $(\sin^2(n\pi a))_{n\in\N}$ $($Definition~$\ref{o1})$ and let $\omega_2$ be defined on $[0,\infty)$ by $\omega_2(t)=t^4.$ Let $\Phi$ and $\Psi$ be two functions such that the quadruplet $(\omega_1,\omega_2,\Phi,\Psi)$ is admissible $($see Definition~$\ref{adm})$ and let $\vH_{\Phi,\Psi}$ be defined by $(\ref{optim1}).$ Then there exists a constant $C=C(a)>0$ such that for any initial data $(u^0,u^1)\in\vD(\vA_3),$ the corresponding solution $u$ of $(\ref{beam})$ verifies
\begin{gather*}
\forall t\ge0,\;
E_3(u(t))\le C\vH_{\Phi,\Psi}\left(\frac{1}{t+1}\right)\|(u^0,u^1)\|_{\vD(\vA_3)}^2,
\end{gather*}
if $\Phi$ and $\Psi$ satisfy the hypothesis~$(\ref{adm3-1})$ of Definition~$\ref{adm}$ and
\begin{gather}
\label{thmdecay0-2-2}
\forall t\ge0,\;
E_3(u(t))\le\frac{C}{(t+1)^\frac{q}{p}}\|(u^0,u^1)\|_{\vD(\vA_3)}^2,
\end{gather}
if for any $t>0,$ $\Phi(t)=C_1t^\frac{1}{p}$ and $\Psi(t)=C_2t^\frac{1}{q}$ for some $p\in[1,\infty),$ $q\in[1,\infty)$ and constants $C_1,C_2>0$ $($case $(\ref{adm3-2})$ of Definition~$\ref{adm}).$
\end{vthm}

\begin{proof*}
We write $u^0(x)=\vsum_{n=0}^\infty a_n\sin(n\pi x)$ and $u^1(x)=\pi^2\vsum_{n=0}^\infty n^2b_n\sin(n\pi x)$ and we consider the solution $v$ of $v_{tt}+v_{xxxx}=0,$ satisfying the same boundary condition and having the same initial data as $u.$ We follow the method as for (\ref{wave}). From Ingham's inequality (Ingham~\cite{MR1545625}) and Lemmas~3.3 and 5.1 of Ammari and Tucsnak~\cite{MR1814271} (see also Lemmas~4.1 and 5.7 of Ammari and Tucsnak~\cite{MR2002f:93104}), we obtain for $T=10,$
\begin{gather*}
\vint_0^{\;T}u_t^2(t,a)dt\ge C(a)\vint_0^{\;T}v_t^2(t,a)dt
\ge C(a)\vsum_{n=0}^\infty n^4(a_n^2+b_n^2)\sin^2(n\pi a).
\end{gather*}
Then we define
\begin{gather*}
E_+(0)=\dfrac{\pi^8}{4}\dsp\vsum_{n=0}^\infty n^8(a_n^2+b_n^2)
=\dfrac{\pi^8}{4}\dsp\vsum_{n=0}^\infty n^4(a_n^2+b_n^2)\omega_2(n), \medskip \\
E(0)=\dfrac{\pi^4}{4}\dsp\vsum_{n=0}^\infty n^4(a_n^2+b_n^2)
\stackrel{{\rm def}}{=}\frac{1}{2}\|(u^0,u^1)\|_{X_3}^2, \medskip \\
E_-(0)=\dfrac{\pi^4}{2}\dsp\vsum_{n=0}^\infty n^4(a_n^2+b_n^2)\omega_1(n).
\end{gather*}
The result follows from the discussion at the beginning of Section~\ref{damping}.
\medskip
\end{proof*}

\noindent
Using Theorem~\ref{thmdecay1-2} and Proposition~\ref{propdefomega} (applied with
$(u_n)_{n\in\N}=(\sin^2(n\pi a))_{n\in\N},$ $p=4$ and $\alpha=0),$ we obtain the following result.

\begin{vthm}
\label{thmdecay2-2}
Let $a\in(0,1)\cap\Q^{\mathrm c}$ and let $\omega_1$ and $\omega_2$ be defined as in Theorem~$\ref{thmdecay1-2}.$ We set
\begin{gather*}
\forall t>0, \; \vphi(t)=\dfrac{\omega_1(t)}{t^4}.
\end{gather*}
Then there exists a constant $C=C(a)>0$ such that for any initial data $(u^0,u^1)\in\vD(\vA_3),$ the solution $u$ of $(\ref{beam})$ satisfies
\begin{gather*}
\forall t\ge0, \; \|(u(t),u_t(t))\|_{V_3\times L^2(0,1)}\le
\dfrac{C}{\left(\vphi^{-1}\left(\frac{1}{t+1}\right)\right)^2} \|(u^0,u^1)\|_{\vD(\vA_3)}.
\end{gather*}
\end{vthm}

\begin{vrmk}
By Theorem~\ref{thmdecay2-2}, we are able to give the explicit decay of the energy for any $a\in(0,1)\cap\Q^{\mathrm c}.$ This completes the lack, since the decay was known for almost every $a\in(0,1)$ (Ammari and Tucsnak~\cite{MR1814271}, Theorem~2.2). In addition, with help of \eqref{thmdecay0-2-2} of Theorem~\ref{thmdecay1-2}, our method allows us to recover the decay of Theorem~2.2 in Ammari and Tucsnak~\cite{MR1814271}.
\medskip
\end{vrmk}

\baselineskip 0.5cm

\noindent
{\large\bf Acknowledgments} \\
The authors are grateful to Professor Enrique Zuazua for stimulating remarks and for pointing out to them several imprecisions in the first version of the paper. The first author would like to thank Professor Enrique Zuazua for having invited him several months in the ``Departamento de Matem\'aticas'' of the ``Universidad Aut\'onoma de Madrid'' and for suggested him this work. Finally, the authors express their gratitude to the referee who made helpful advices, which permitted to clarify and improve some points of this paper. He also made some useful suggestions to extend the decay results of Sections~\ref{damping} and \ref{damping2} to higher dimensions.

\baselineskip 0cm

\bibliographystyle{abbrv}
\bibliography{Paper5}
\addcontentsline{toc}{section}{References}

\end{document}